\documentstyle[12pt,amscd]{amsart}

\font\chq = cmr7 
\font\rechq = cmr5 scaled \magstep 0

\font\tenmsb=msbm10 scaled \magstep 1
\font\sevenmsb=msbm7 scaled \magstep 1
\newfam\msbfam
\textfont\msbfam=\tenmsb
\scriptfont\msbfam=\sevenmsb
\def\Bbb{\fam\msbfam\tenmsb}

\font\tenfrak=eufm10 scaled \magstep1
\font\sevenfrak=eufm7 scaled \magstep1
\newfam\frakfam
\textfont\frakfam=\tenfrak
\scriptfont\frakfam\sevenfrak
\def\frak{\fam\frakfam\tenfrak}

\numberwithin{equation}{subsection}
\theoremstyle{plain}
	\begingroup
	\newtheorem{thm}[equation]{Theorem}
	\newtheorem{lem}[equation]{Lemma}
	\newtheorem{prop}[equation]{Proposition}
	\newtheorem{cor}[equation]{Corollary}
	\newtheorem{defn}[equation]{Definition}
	\newtheorem{exmp}[equation]{Example}
\endgroup
\theoremstyle{remark}
	\newtheorem{rem}[equation]{Remark}

\newsymbol\leftrightarrows 131C
\newsymbol\rightleftarrows 131D
\newsymbol\rightrightarrows 1313

\newcommand{\so}{\Rightarrow}

\newcommand{\sii}{\Leftrightarrow}
\newsymbol\twoheadrightarrow 1310
\newcommand{\sobre}{\twoheadrightarrow}             
\newcommand{\gania}{\leftharpoonup}                 
\newcommand{\ganda}{\rightharpoonup}                
\newcommand{\ganid}{\leftharpoondown}               
\def\muyrarrow{\DOTSB\relbar\joinrel\relbar\joinrel\relbar\joinrel\rightarrow}
\def\muymuyrarrow{\DOTSB\relbar\joinrel\relbar\joinrel\relbar\joinrel\relbar\joinrel\rightarrow}
\def\muymuymuyrarrow{\DOTSB\relbar\joinrel\relbar\joinrel\relbar\joinrel\relbar\joinrel\relbar\joinrel\rightarrow}
\newcommand{\flad}[1]{\displaystyle
		\mathop{\rightarrow}^{#1}}          
\newcommand{\fllad}[1]{\displaystyle
		\mathop{\longrightarrow}^{#1}}      
\newcommand{\flllad}[1]{\displaystyle
		\mathop{\muyrarrow}^{#1}}           
\newcommand{\fllllad}[1]{\displaystyle
		\mathop{\muymuyrarrow}^{#1}}        
\newcommand{\flllllad}[1]{\displaystyle
		\mathop{\muymuymuyrarrow}^{#1}}     

\newcommand{\ev}{\operatorname{ev}}

\newcommand{\Hom}{\operatorname{Hom}}
\newcommand{\End}{\operatorname{End}}
\newcommand{\Aut}{\operatorname{Aut}}

\newcommand{\Ad}{\operatorname{Ad}}
\newcommand{\equ}{\operatorname{Eq}}
\newcommand{\id}{\operatorname{id}}

\newcommand{\LKer}{\operatorname{LKer}}

\newcommand{\Ind}{\operatorname{Ind}}
\newcommand{\sg}{\operatorname{sg}}
\newcommand{\gra}{\operatorname{gr}}
\newcommand{\car}{\operatorname{char}}
\newcommand{\Ss}{{\cal S}}
\newcommand{\St}{{\cal S}^{-1}}
\newcommand{\TTT}{{\bf T}}                          
\newcommand{\SSS}{{\bf S}}                          

\newcommand{\ZZ}{{\Bbb Z}}
\newcommand{\NN}{{\Bbb N}}
\renewcommand{\SS}{{\Bbb S}}
\newcommand{\BB}{{\Bbb B}}
\newcommand{\DD}{{\Bbb D}}
\renewcommand{\k}{{\bold k}}
\newcommand{\h}{{\bold h}}
\newcommand{\f}{{\bold f}}
\newcommand{\bo}{{\frak b}}
\newcommand{\g}{{\frak g}}
\newcommand{\n}{{\frak n}}
\renewcommand{\u}{{\frak u}}
\newcommand{\toba}[1]{{\frak t}(#1)}		    
\newcommand{\tobag}[2]{{\frak t}^{#1}(#2)}	    
\newcommand{\U}{{\bold U}}
\newcommand{\ide}{{\cal I}}                         
\newcommand{\cate}{{\cal C}}                        
\newcommand{\cateb}{\overline{{\cal C}}}            
\newcommand{\modui}[1]{{}_{#1}\!{\cal M}}           
\newcommand{\modud}[1]{{\cal M}_{#1}}               
\newcommand{\comoi}[1]{{}^{#1}\!{\cal M}}           
\newcommand{\comod}[1]{{\cal M}^{#1}}               
\newcommand{\Modui}[1]{{}_{#1}^\infty\!{\cal M}}    
\newcommand{\Modud}[1]{{\cal M}_{#1}^\infty}        
\newcommand{\hopmoi}[1]{{}_{#1}^{#1}{\cal M}}       
\newcommand{\hopf}[1]{{\cal H}f(#1)}                
\newcommand{\uni}{{\bold1}}                         
\newcommand{\yd}{{}^H_H{\cal YD}}                   
\newcommand{\ydh}{{\cal YD}^H_H}                    
\newcommand{\yds}[1]{{}^{#1}_{#1}{\cal YD}}         
\newcommand{\ydb}{{}^H_H\overline{{\cal YD}}}       
\newcommand{\ydi}{{}^H_H{\cal YD}_\infty}           
\newcommand{\hsyd}{{}^{H^*}_{H^*}{\cal YD}}         
\newcommand{\ydhs}{{\cal YD}^{H^*}_{H^*}}           
\newcommand{\hbyd}{{}^{H^{\hbox{\rechq bop}}}_{H^{\hbox{\rechq bop}}}{\cal YD}}     
\newcommand{\ydhb}{{\cal YD}^{H^{\hbox{\rechq bop}}}_{H^{\hbox{\rechq bop}}}}       
\newcommand{\ydhc}{{\cal YD}^{H^{\hbox{\rechq cop}}}_{H^{\hbox{\rechq cop}}}}       
\newcommand{\hsbyd}{{}^{H^{*\hbox{\rechq bop}}}_{H^{*\hbox{\rechq bop}}}{\cal YD}}  
\newcommand{\ydhsb}{{\cal YD}^{H^{*\hbox{\rechq bop}}}_{H^{*\hbox{\rechq bop}}}}    
\newcommand{\hobim}{{}^H_H{\cal M}^H_H}             
\newcommand{\doble}[1]{#1\!\bowtie #1^{*\hbox{\chq op}}}       
\newcommand{\dob}[1]{{\cal D}(#1)}                               
\newcommand{\orb}{{\cal O}}                         
\newcommand{\conor}{{\cal Q}}                       
\newcommand{\rid}[1]{\widehat{#1}}                  
\newcommand{\prim}[1]{{\cal P}(#1)}                 
\newcommand{\hop}{H^{\hbox{\chq op}}}               
\newcommand{\hcop}{H^{\hbox{\chq cop}}}             
\newcommand{\hbop}{H^{\hbox{\chq bop}}}             
\newcommand{\hsop}{H^{*\hbox{\chq op}}}             
\newcommand{\aop}{A^{\hbox{\rechq op}}}             
\newcommand{\ccop}{C^{\hbox{\rechq cop}}}           
\newcommand{\rop}{R^{\hbox{\rechq op}}}             
\newcommand{\rcop}{R^{\hbox{\rechq cop}}}           
\newcommand{\rbop}{R^{\hbox{\rechq bop}}}           
\newcommand{\rsbop}{R^{*\hbox{\rechq bop}}}         
\newcommand{\invi}[2]{{}^{#2}#1}                    
\newcommand{\invd}[2]{#1^{#2}}                      
\newcommand{\coini}[2]{{}^{\hbox{\chq co}#2}#1}     
\newcommand{\coind}[2]{#1^{\hbox{\chq co}#2}}       
\newcommand{\intgl}[1]{{\cal I}_{\ell}(#1)}         
\newcommand{\intgr}[1]{{\cal I}_r(#1)}              
\newcommand{\I}{{\cal I}}                           
\newcommand{\dai}[1]{{}^*\!#1}                      
\newcommand{\dad}[1]{#1^*}                          
\newcommand{\bdh}{{}_{\alpha}h}                     
\newcommand{\dbh}{{}^{\alpha}h}                     
\newcommand{\bdm}{{}_{\alpha}m}                     
\newcommand{\dbm}{{}^{\alpha}m}                     
\newcommand{\bdr}{{}_{\alpha}r}                     
\newcommand{\dbr}{{}^{\alpha}r}                     
\newcommand{\bdrb}{{}_{\beta}r}                     
\newcommand{\dbrb}{{}^{\beta}r}                     
\newcommand{\agg}[2]{#1(#2)}                        

\newcommand{\alai}[2]{{}^{#2}#1}                    
\newcommand{\estri}{{\cal R}}                       
\newcommand{\cou}{\varepsilon}                      
\newcommand{\com}[2]{#1_{(#2)}}                     
\newcommand{\comm}[3]{#1_{(#2)(#3)}}                
\newcommand{\commm}[4]{#1_{(#2)(#3)(#4)}}           
\newcommand{\commn}[4]{#1_{(#2)(#3)(-#4)}}          
\newcommand{\con}[2]{#1_{(-#2)}}                    
\newcommand{\comn}[3]{#1_{(#2)(-#3)}}               
\newcommand{\mulop}{m^{\chq\hbox{op}}}              
\newcommand{\coporop}{\Delta^{\chq\hbox{op}}}       
\newcommand{\nftp}[3]{#1_{#2}\otimes\ldots\otimes #1_{#3}} 

\newcommand{\noi}{\noindent}
\newcommand{\gi}{\langle}                           
\newcommand{\gd}{\rangle}                           
\newcommand{\oti}{{\scriptstyle\otimes}}            
\newcommand{\coti}{\Box}                            
\newcommand{\ba}[1]{\overline{#1}}                  
\newcommand{\gp}[1]{\!<\!\!#1\!>}                   

\begin{document}

\title[Braided Hopf algebras]
	{Braided Hopf algebras over non abelian finite groups}
\author{Nicol\'as Andruskiewitsch}
\address{Nicol\'as Andruskiewitsch \\ FaMAF -- C\'ordoba}  \email{andrus@@mate.uncor.edu} \thanks{This work was partially supported by CONICOR, CONICET, SECyT-UNC and TWAS}

\author{Mat\'\i as Gra\~na}
\address{Mat\'\i as Gra\~na \\ FCEyN -- Buenos Aires} \email{matiasg@@dm.uba.ar}
\date{May 20, 1998. {\tt 9802074}}
\maketitle

\begin{abstract}In the last years a new theory of Hopf algebras has begun
to be developed: that of Hopf algebras in braided categories, or, briefly,
braided Hopf algebras. This is a survey of general aspects of the theory
with emphasis in $\yd$, the Yetter--Drinfeld category over $H$,
where $H$ is the group algebra of a non abelian finite group $\Gamma$.
We discuss a special class of braided graded Hopf algebras from
different points of view following Lusztig, Nichols and Schauenburg.
We present some finite dimensional examples arising in an unpublished
work by Milinski and Schneider.
\\[.1in]
{\sc Sinopsis.} En los \'ultimos a\~nos comenz\'o a ser desarrollada una nueva
teor\'\i a de \'algebras de Hopf en categor\'\i as trenzadas, o brevemente,
\'algebras de Hopf trenzadas. Presentamos aqu\'{\i} aspectos generales de la
teor\'\i a con \'enfasis en $\yd$, la categor\'\i a de Yetter--Drinfeld
sobre $H$, donde $H$ es el \'algebra de grupo de un grupo finito no abeliano
$\Gamma$.
Discutimos una clase especial de \'algebras de Hopf trenzadas graduadas
desde diferentes puntos de vista, siguiendo a Lusztig, Nichols y
Schauenburg. Presentamos algunos ejemplos de dimensi\'on finita que
aparecen en un trabajo in\'edito de Milinski y Schneider.
\end{abstract}

\addtocounter{section}{-1}

\section{Introduction and notations}
\subsection{Introduction} \ \\
The idea of considering Hopf algebras in braided categories (categories
with a tensor product which is associative and ``commutative") goes back
to Milnor--Moore \cite{milnor-moore} and Mac Lane \cite{maclane:63}.
Hopf superalgebras, or $\ZZ/2$-graded Hopf algebras, were intensively
studied in the work of Kac, Kostant, Berezin and others.
In this case, the braiding $c$ is symmetric: $c^2=\id$. With the advent
of quantum groups, it became clear that braided (non symmetric) categories
have a r\^ole to play in several parts of algebra. This point of view was
pioneered by Manin \cite[p. 81]{manin:88}, Majid \cite{majid:95}
(see also \cite{gurevich:91}) and widely developed since then.
One of its main applications is Lusztig's presentation of quantized enveloping
algebras \cite{lusztig:93}.

Our motivation to study braided Hopf algebras is the so-called bosonization
(or biproduct) construction, due to Radford \cite{radford:85} and
interpreted in the terms of braided categories by Majid \cite{majid:94a}.
More precisely, we are interested in a specific type of braided Hopf
algebras. To explain the reason we recall a general principle from
\cite{andrus-schneider:97}:

Let $K$ be a Hopf algebra with coradical filtration
$$K_0\subset K_1\subset\ldots.$$
If the coradical $K_0$ is a Hopf subalgebra (this happens for instance
if $K$ is pointed, in which case $H=K_0$ is a group algebra)
then the associated graded space
$$\gra K=\bigoplus_{n\ge 0}{K_n}/{K_{n-1}}\qquad (K_{-1}=0)$$
has a graded Hopf algebra structure inherited from that of $K$.
Moreover, since the inclusion \linebreak
${K_0\hookrightarrow\gra K}$ has a
retraction $\gra K\to K_0$ of Hopf algebras, the inverse process
to the bosonization construction makes the algebra of coinvariants
$R=\coind{(\gra K)}{K_0}$ into a graded Hopf algebra in $\yds {K_0}$
with trivial coradical. Conversely, if $H$ is a group algebra,
let $R$ be a graded Hopf algebra in $\yd$ with trivial coradical
(i.e., $R_0=\k 1$). Then the bosonization $R\# H$ is a pointed graded
Hopf algebra with coradical isomorphic to $H$.
It is then reasonable to expect that information one can give about graded
Hopf algebras in $\yd$ can be translated to information about pointed Hopf
algebras.

We say that a graded Hopf algebra $R=\bigoplus_{i\ge 0}\agg Ri$ in $\yd$
is a TOBA if $\agg R0$ is the base field (and then the coradical is trivial by
\cite[11.1.1]{sweedler:69}), the space of primitive elements
is exactly $\agg R1$ and this space generates $R$ (as an algebra).
It is then proved that $R$ is a TOBA iff $R\# H$ is a Hopf algebra of
type one, in the sense of Nichols \cite{nichols:78}.
An important example of TOBA is the quantum analog of the enveloping algebra
of the nilpotent part of a Borel algebra, see \cite{lusztig:93},
\cite{schauenburg:96}, \cite{rosso:95}, \cite{rosso:92}.

The article is organized as follows:

In section 1 we define and give examples of braided categories, Hopf
algebras in braided categories and review the bosonization construction.

In section 2 we give duals and opposite algebras of braided Hopf
algebras (a deep treatment of the subject can be found in
\cite{majid:95}). For finite Hopf algebras we define the space of
integrals, and prove (following Takeuchi \cite{takeuchi:97}) that it is an
invertible object in the category. This allows to state ``braided"
versions of several useful results concerning finite dimensional Hopf
algebras (e.g. the bijectivity of the antipode).

In section 3 we concentrate on braided Hopf algebras in $\yd$ where
$H=\k\Gamma$, the group algebra of a finite group $\Gamma$.
We show that a TOBA $R$ is determined, up to isomorphism, by the space
of primitive elements $\prim R=\agg R1$. Moreover, given a Yetter-Drinfeld
module $V$, we present three different constructions of a TOBA $\toba V$
such that its space of primitive elements is isomorphic to $V$. The first
two constructions use quantum shuffles and universal properties, and are
essentially contained in \cite{nichols:78}, \cite{schauenburg:93},
\cite{rosso:95}, \cite{rozanski:96}, \cite{woronowicz:89}.
The third construction, by means of a bilinear form,
seems to be new. It is however inspired by \cite{lusztig:93},
\cite{schauenburg:93}, \cite{rosso:95}, \cite{muller:98}.
We finally discuss some explicit examples from \cite{milinski-schneider}
for $\Gamma$ a symmetric or a dihedral group.

We thank H.-J. Schneider for valuable conversations
during the preparation of this article.

\subsection{Notations} \ \\
We shall work over a field $\k$. Sometimes we impose some hypothesis to
the field, but in most of the article it may be any field.
Tensor products and Homs are taken over $\k$ when not specified.
We use the letters $H,K$ for Hopf algebras over $\k$, and the letter $R$
for braided Hopf algebras.

Given a Hopf algebra $H$, we use subindices $H_0\subset H_1\subset\ldots$
to indicate the coradical filtration of $H$ (see \cite{sweedler:69}).
In order to avoid confusion with this notation, a graded algebra shall
be denoted by $H=\oplus_i\agg Hi$.

Given an algebra $A$, we denote by $\modui A$ the category of
finite dimensional left $A$-modules and by $\Modui A$ the category
of all left $A$-modules; {\it ditto} for the categories of right
modules $\modud A$ and $\Modud A$.

Given a coalgebra $C$, we denote by $\comoi C$ the category of
left $C$-comodules. The same for $\comod C$.

Given a Hopf algebra $H$, we denote by $\hobim$ the category of Hopf
bimodules over $H$, by $\yd$ the Yetter--Drinfeld category over $H$
(of finite dimensional left YD modules) and by $\ydi$ the category of
all left YD modules. See \ref{defyd} for the precise definitions.

Given a Hopf algebra $H$ with bijective antipode, we denote by $\hop$ the
Hopf algebra with opposite multiplication, $\hcop$ the Hopf algebra with
opposite comultiplication, and $\hbop$ the Hopf algebra
$(\hop)^{\hbox{\chq cop}}$.

We use for coalgebras Sweedler notation without summation symbol:
$\Delta(h)=\com h1\otimes\com h2$,
and the same for comodules: $\delta(m)=\con m1\otimes\com m0$ for left
comodules and $\delta(m)=\com m0\otimes\com m1$ for right comodules.

If $M$ is a $\k$-vector space, $m\in M$ and $f\in\dad M$, we use either
$f(m)$, $\gi f,m\gd$ or $\gi m,f\gd$ to denote the evaluation map.

\section{Definitions and examples}
\subsection{Braided categories. Abelian braided categories} \ \\
\label{predef}
We recall in this section the definitions of monoidal and braided
categories. See \cite{joyal-street} for a detailed treatment of the subject.

\begin{defn}
A {\em monoidal category} is a category $\cate$ together with a functor
$\otimes:\cate\times\cate\to\cate$ (called {\em tensor product}), an object
$\uni$ of $\cate$ (called {\em unit}) and natural isomorphisms
\begin{alignat*}{2}
&a_{U,V,W}:(U\otimes V)\otimes W\to U\otimes(V\otimes W) \qquad
        && \mbox{(associativity constraint)}, \\
&r_V:V\otimes \uni\to V,\qquad l_V:\uni\otimes V\to V \qquad
        && \mbox{(unit constraints)},
\end{alignat*}
subject to the following conditions:
\begin{multline*}
((U\otimes V)\otimes W)\otimes X\to (U\otimes V)\otimes(W\otimes X)\to
U\otimes(V\otimes(W\otimes X)) \\
=((U\otimes V)\otimes W)\otimes X \to (U\otimes(V\otimes W))\otimes X\to
U\otimes((V\otimes W)\otimes X)\to U\otimes(V\otimes(W\otimes X)),
\end{multline*}
$$\id=V\otimes W\to(V\otimes \uni)\otimes W\to V\otimes (\uni\otimes W)
\to V\otimes W.$$
\end{defn}

We shall assume in what follows that the associativity constraint is
the identity morphism. This is possible thanks to \cite{maclane:71},
where the author proves that any monoidal category can be embedded in
another monoidal category in which this is true.

\begin{defn}
A {\em braided monoidal category} (or briefly {\em braided category})
is a monoidal category $\cate$ together with a natural isomorphism
$$c=c_{M,N}:M\otimes N\to N\otimes M$$
(called {\em braiding}) subject to the conditions
\begin{align}
c_{M,N\otimes P}&=(\id_N\otimes c_{M,P})\circ
        (c_{M,N}\otimes\id_P), \label{eqbt1} \\
c_{M\otimes N,P}&=(c_{M,P}\otimes\id_N)\circ
        (\id_M\otimes c_{N,P}). \label{eqbt2}
\end{align}
The category $\cate$ is called {\em symmetric} if $c^2=\id$, i.e.,
for all $M,N$ in $\cate$, $c_{N,M}c_{M,N}=\id_{M\otimes N}$.
\end{defn}

There are some identities that can be proved from the axioms (and
hence hold in any braided category). One of these identities is
$lc=r$, that is,
$$(M \otimes \uni @>c_{M,\uni}>> \uni\otimes M @>l_M>> M) =
(M \otimes \uni @>r_M>> M).$$
Analogously $rc=l$.

\begin{rem} \label{rmfunmon}
To define when two monoidal (braided) categories are equivalent, it is
necessary to know what a functor between monoidal (braided) categories
is. Let $\cate$ and $\cate'$ be monoidal categories. A functor between
them is a pair $(F,\eta)$, where $F:\cate\to\cate'$ is a functor
and $\eta$ is a natural isomorphism
$\eta:\otimes\circ F^2\to F\circ\otimes$ (that is,
$\eta_{M,N}:FM\otimes FN\to F(M\otimes N)$) subject to the conditions
\begin{equation} \label{funmon}
  \begin{CD}
  FM\otimes FN\otimes FP @>>> F(M\otimes N)\otimes FP \\
  @VVV @VVV \\
  FM\otimes F(N\otimes P) @>>> F(M\otimes N\otimes P)
  \end{CD}
        \quad\hbox{must commute,}
\end{equation}
\begin{align*}
F(\uni_\cate) &= \uni_{\cate'}, \\
\uni\otimes FM \fllad{l_{FM}} FM &= \uni\otimes FM \flad{\eta}
        F(\uni\otimes M) @>F(l_M)>> FM, \\
FM\otimes \uni \fllad{r_{FM}} FM &= FM\otimes\uni \flad{\eta}
        F(M\otimes\uni) \fllad{F(r_M)} FM.
\end{align*}
Observe that when the associativity constraint is not the identity,
then \eqref{funmon} must be suitably modified. For braided categories
the following diagram must also commute.
$$\begin{CD}
FM\otimes FN @>c_{FM,FN}>> FN\otimes FM \\
@V\eta VV @V\eta VV \\
F(M\otimes N) @>F(c_{M,N})>> F(N\otimes M).
\end{CD}$$
\end{rem}

\bigskip

Given a braided category $\cate$, we shall denote by $\cateb$ the braided
category whose objects and morphisms are those of $\cate$ but whose braiding
is the inverse of that of $\cate$. The axioms \eqref{eqbt1} and \eqref{eqbt2}
are automatically verified for this category.

\medskip

An important source of examples of braided categories is given by the
quasitriangular bialgebras. Let $H$ be a bialgebra. An element
$\estri\in H\otimes H$ is called
a {\em triangular structure} for $H$ if it is invertible (with respect to
the usual product of $H\otimes H$) and verifies
\begin{align}
\forall h\in H,\ \coporop(h)&=\estri\Delta(h)\estri^{-1}, \label{eqet1} \\
(\id\otimes\Delta)(\estri)&=\estri^{13}\estri^{12}, \label{eqet2} \\
(\Delta\otimes\id)(\estri)&=\estri^{13}\estri^{23}. \label{eqet3}
\end{align}
In this case, the pair $(H,\estri)$ is called a {\em quasitriangular
bialgebra} (QT bialgebra). If $\tau(\estri)=\estri^{-1}$
($\tau$ is the usual flip), then $(H,\estri)$ is called {\em triangular}.
For $(H,\estri)$ a QT bialgebra, the category of left (right) $H$-modules
$\Modui H$ ($\Modud H$) and the category of finite dimension left (right)
$H$-modules $\modui H$ ($\modud H$) are braided, where
\begin{align*}
c_{M,N}(m\otimes n)=\estri_2n\otimes\estri_1m\qquad &\mbox{for left modules}, \\
c_{M,N}(m\otimes n)=n\estri_1\otimes m\estri_2\qquad &\mbox{for right modules}.
\end{align*}

Equation \eqref{eqet1} is equivalent to $c$ being a morphism of $H$-modules.
Equations \eqref{eqet2} and \eqref{eqet3} are respectively equivalent to
\eqref{eqbt1} and \eqref{eqbt2} in the case of left modules, and to
\eqref{eqbt2} and \eqref{eqbt1} in the case of right modules.
These categories are symmetric if $(H,\estri)$ is triangular.

The notion of QT bialgebra can be dualized to that of
{\em co-quasitriangular bialgebras} (or briefly CQT bialgebras), for which
the category of left (right) comodules is braided.
Both notions can be generalized to quasi-bialgebras, to get QT
{\em quasi-bialgebras} (CQT quasi-bialgebras) for which the categories
of left (right) modules (left (right) comodules) are also braided
(see \cite{drinfeld:90b}).
In these cases the associativity constraint is no longer the usual
associativity for vector spaces, and the verification of the axioms
becomes more tedious.

\begin{defn} \label{duales}
A monoidal category $\cate$ is called {\em rigid} if every object has
a left and right dual in it. That is, for every object $M$ of $\cate$
there exist $\dai M$ and $\dad M$ objects of $\cate$ and natural morphisms
\begin{align*}
br_M &: \uni\to M\otimes \dad M, \\
bl_M &: \uni\to \dai M\otimes M, \\
dr_M &: \dad M\otimes M\to \uni, \\
dl_M &: M\otimes \dai M\to \uni,
\end{align*}
subject to the conditions
\begin{equation} \label{eqdual} \begin{split}
\id=M \flllad{\id\otimes bl} M\otimes\dai M\otimes M
        \flllad{dl\otimes\id} \uni\otimes M \flad{l} M, \\
\id=M \flllad{br\otimes\id} M\otimes\dad M\otimes M
        \flllad{\id\otimes dr} M\otimes \uni \flad{r} M.
\end{split} \end{equation}
\end{defn}

\begin{rem} \label{buedef*}
The conditions on $\dad M$ and $\dai M$ determine them up to isomorphism.
We shall use in what follows the word ``rigid" for a category in which
the correspondences $M\mapsto\dad M$ and $M\mapsto\dai M$ are given
by functors, which is true in the usual cases.
\end{rem}

\begin{defn} \label{defbn}
It is well known that the symmetric group in $n$ elements $\SS_n\ (n\ge 2)$,
can be presented by elementary transpositions $\tau_i=(i,i+1),\ (1\le i<n)$
subject to the relations
\begin{alignat*}{2}
\tau_i\tau_j &= \tau_j\tau_i\quad && \mbox{if }|i-j|>1, \\
\tau_i\tau_j\tau_i &= \tau_j\tau_i\tau_j\quad && \mbox{if }|i-j|=1, \\
\tau_i^2 &= 1\quad && \forall i.
\end{alignat*}
If we drop the last set of relations, we get the Artin {\em braid group}.
To be precise, we define $\BB_n,\ (n\ge 2)$ to be the group with generators
$\sigma_i,\ (1\le i<n)$ subject to the relations
\begin{alignat*}{2}
\sigma_i\sigma_j &=\sigma_j\sigma_i\quad && \mbox{if }|i-j|>1, \\
\sigma_i\sigma_j\sigma_i &= \sigma_j\sigma_i\sigma_j\quad && \mbox{if }|i-j|=1.
\end{alignat*}
This is an infinite group, and there is a projection map from $\BB_n$ to
$\SS_n$ given by $\sigma_i\mapsto\tau_i$.
\end{defn}

\medskip

\begin{rem} \label{bnactua}
The group $\SS_n$ acts naturally on $n$-fold tensor products in the category
of vector spaces, or more generally, in the category of representations of
a cocommutative Hopf algebra. In both cases the category is symmetric. When
this does not happen, the group $\SS_n$ has to be replaced by $\BB_n$, as we
now explain. Let $\cate$ be a braided category and $M$ an object of $\cate$.
Then $\BB_n$ acts on
$\underbrace{M\otimes\dots\otimes M}_{n\hbox{\chq{ times}}}$
via $$\sigma_i\mapsto \underbrace{\id\otimes\dots\otimes\id}_{i-1}\otimes c
\otimes\underbrace{\id\otimes\dots\otimes\id}_{n-i-1}.$$
This useful observation allows to translate several statements into
drawings, and in fact many authors do use drawings to prove certain
equalities. The axioms above can be viewed as rules to pass from one
configuration to another.
\end{rem}

\bigskip

Our main example of braided (rigid) category is the Yetter--Drinfeld
category over a Hopf algebra:

\begin{defn} \label{defyd}
Let $H$ be a Hopf algebra over $\k$ with bijective antipode.
We shall denote by $\yd$ the category of finite left
{\em Yetter--Drinfeld modules} over $H$. That is, $M$ is an object in
$\yd$ if $M$ is a left $H$-module, a left $H$-comodule, has finite dimension
over $\k$ and
$$\con{(hm)}1 \otimes\com{(hm)}0=\com h1\con m1\Ss\com h3\otimes
\com h2\com m0,\qquad\forall h\in H,\ m\in M.$$
$\yd$ is a monoidal category with the usual tensor product over $\k$,
where $\uni=\k$ and associativity and unit constraints are the usual ones
for vector spaces and where, for $M,N\in\yd$, $M\otimes N$
has the diagonal module and comodule structures given by
$$h(m\otimes n)=\com h1m\otimes\com h2n,\quad
\con {(m\otimes n)}1 \otimes\com{(m\otimes n)}0=
\con m1\con n1\otimes \com m0\otimes\com n0.$$
It is also a braided category, where the braiding is given by
$$c=c_{M,N}:M\otimes N\to N\otimes M,\quad
c(m\otimes n)=\con m1 n\otimes\com m0.$$
It is immediate to see that $c$ is an isomorphism, with inverse
$$(c_{N,M})^{-1}=(c^{-1})_{M,N}\quad
        m\otimes n\mapsto\com n0\otimes\St(\con n1)m.$$
As with any braided category, we shall denote by $\ydb$ the same category
as $\yd$ but with the inverse braiding.
$\ydb$ is a Yetter--Drinfeld category (see \ref{cateigua}).
We prove rigidity of these categories in \ref{ydesri}.
We denote by $\ydi$ the category of all (non necessarily finite dimensional)
Yetter--Drinfeld modules over $H$. This is a braided category (with the
braiding given by the same formula as in $\yd$) but it is not rigid.

\medskip

If $M$ is an object in $\yd$, we shall denote by $\Im M$ the object of $\yd$
with the same underlying vector space but with the structure given by
$$h\ganda m=\Ss^2(h)m,\ \ \delta_{\Im M}(m)=\Ss^{-2}\con m1\otimes\com m0.$$
\end{defn}

\begin{thm}[Majid]
\label{doble}
Let $H$ be a finite dimensional Hopf algebra.
Let $\dob H=\doble H$ be the Drinfeld double of $H$, defined by
$\dob H=H\otimes H^*$ as a coalgebra, with multiplication and antipode
given by (we denote here $fg=m(f\otimes g)$ in $H^*$ rather than in $\hsop$)
\begin{align*}
(h\bowtie f)(h'\bowtie f')&=\gi\com f1,\com {h'}1\gd\gi\com f3,\Ss\com {h'}3\gd
(h\com {h'}2\bowtie f'\com f2), \\
\Ss_{\dob H}(h\bowtie f)&=(1\bowtie\St f)(\Ss h\bowtie\cou)
=(\Ss\com h2\bowtie\St\com f2)\gi\com f1,\Ss\com h1\gd\gi\com f3,\com h3\gd.
\end{align*}
We observe that $H$ and $\hsop$ are Hopf subalgebras of $\dob H$.
We have that $\dob H$ is a QT Hopf algebra, and the category $\modui{\dob H}$
of finite dimensional left $\dob H$-modules is equivalent, as braided
category, to $\yd$.
\end{thm}

\begin{pf} See \cite{montgomery:93}. \end{pf}

There is another category which appears naturally in the framework of Hopf
algebras which is equivalent to $\yd$ and $\modui{\dob H}$, namely $\hobim$.
This is the category whose objects are $H$-bimodules
and $H$-bicomodules, such that the structure morphisms $H\otimes M\to M$
and $M\otimes H\to M$ are bicomodule morphisms, taking in $M\otimes H$ and
$H\otimes M$ the codiagonal structure (equivalently, the structure
morphisms $M\to M\otimes H$ and $M\to H\otimes M$ are bimodule morphisms
taking in $M\otimes H$ and $H\otimes M$ the diagonal structure). We take
in this category tensor products over $H$ (alternatively, we can take
the monoidal structure given by cotensor products over $H$). This category
has a braiding, namely
$$c_{M,N}(m\otimes n)=
        \con m2\com n0\Ss(\com n1)\Ss(\con m1)\otimes\com m0\com n2.$$

The following result was independently found by Schauenburg and
the first author.
\begin{prop} \label{yd-hobim}
The category $\hobim$ is equivalent as braided category to $\yd$.
(Alternatively, the category with the monoidal structure given by
cotensor products is also equivalent).
\end{prop}

\begin{pf}(Sketch, see \cite[Satz 1.3.5]{schauenburg:93},
\cite{schauenburg:94} or \cite[Appendix]{andrus-devoto:95} for the details).
Let $M$ be in $\hobim$.
By \cite[Th. 4.1.1]{sweedler:69}, $M\simeq V\otimes H$ as a
right module and right comodule, where
$V=M^{\hbox{\chq co }H}$ and the right module and comodule structures of
$V\otimes H$ are those of $H$. Let us identify $V$ with $V\otimes 1$.
We take the structure of $V$ in $\yd$ by
$$h\ganda v=\com h1v\Ss\com h2,\quad\delta(v)=\delta_l(v)
        \quad\mbox{($\delta_l:M\to H\otimes M$ is the structure morphism).}$$

Conversely, if $V$ is a Yetter--Drinfeld module over $H$, then
$V\otimes H$ is an object in $\hobim$ via
$$h(v\otimes g)=\com h1\ganda v\otimes\com h2g,\quad
        \delta(v\otimes g)=\con v1\com g1\otimes(\com v0\otimes\com g2),$$
and $H$ acts and coacts on the right only over $H$.
\end{pf}

\subsection{Hopf algebras in braided categories (braided Hopf algebras)} \ \\
Monoidal categories are the natural context to define algebras and coalgebras.
If $\cate$ is a monoidal category, we define an algebra in $\cate$ to be
a pair $(A,m)$, where $A$ is an object of $\cate$, $m:A\otimes A\to A$ is
a morphism in $\cate$ and there exists a morphism $u:\uni\to A$ such that
\begin{align*}
m\circ (m\otimes\id) &= m\circ (\id\otimes m):A\otimes A\otimes A\to A, \\
m\circ (u\otimes\id)\circ l_A^{-1} &= \id =
        m\circ (\id\otimes u)\circ r_A^{-1}:A\to A.
\end{align*}
We define dually a coalgebra in $\cate$ to be a pair $(C,\Delta)$, where
$C$ is an object of $\cate$, $\Delta:C\to C\otimes C$ is a morphism in
$\cate$ and there exists a morphism $\cou:C\to\uni$ such that
\begin{align*}
(\Delta\otimes\id)\circ\Delta &= (\id\otimes\Delta)\circ\Delta
        :C\to C\otimes C\otimes C, \\
l_C\circ (\cou\otimes\id)\circ \Delta &= \id =
        r_C\circ (\id\otimes\cou)\circ \Delta:C\to C.
\end{align*}

\medskip

In turn, braided categories are the natural context to define bialgebras
and Hopf algebras.
\begin{defn}
Let $\cate$ be a braided category. A {\em bialgebra in $\cate$} is a triple
$(R,m,\Delta)$, where $R$ is an object in $\cate$ and there exist morphisms
$u:\uni\to R$ and $\cou:R\to\uni$ in such a way that
$(R,u,m)$ is an algebra in $\cate$, $(R,\cou,\Delta)$ is a coalgebra
in $\cate$, $\cou$ is an algebra morphism,
and the usual compatibility between $m$ and $\Delta$ is replaced by
$$\Delta m=(m\otimes m)\circ
        (\id_R\otimes c_{R,R}\otimes\id_R)\circ(\Delta\otimes\Delta).$$
If moreover there exists a morphism $\Ss:R\to R$ which is the inverse
of the identity in the monoid $\Hom_{\cate}(R,R)$ with the convolution
product, then we say that $R$ is a {\em Hopf algebra} in $\cate$ and call
$\Ss$ the {\em antipode} of $R$. We recall that this last definition
can be stated in other words as
$$m(\Ss\otimes\id)\Delta=u\cou=m(\id\otimes\Ss)\Delta:R\to R.$$
\end{defn}

As in the classical case, the compatibility between the algebra and coalgebra
structure can be alternatively stated saying that $m$ is a morphism of
coalgebras, or that $\Delta$ is a morphism of algebras, with the only
difference that $R\otimes R$ is considered as an algebra with the product 
$${m_{R\otimes R}=}{(m_R\otimes m_R)}\circ
{(\id_R\otimes c_{R,R}\otimes \id_R)}$$
or as a coalgebra with the coproduct
$${\Delta_{R\otimes R}=}{(\id_R\otimes c_{R,R}\otimes\id_R)}\circ
{(\Delta_R\otimes\Delta_R)}.$$

As in the classical case, the antipode of a braided Hopf algebra
twists multiplications and comultiplications. One should be careful
to distinguish between $c$ and $c^{-1}$.
The precise equalities are the following ones.
\begin{lem} \label{sdavuel}
Let $R$ be a Hopf algebra in a braided category. Let us denote 
$m=m_R$, $\Delta=\Delta_R$, $c=c_{R,R}$. Then 
$$\Ss_R m=m(\Ss_R\otimes\Ss_R)c, \qquad
        \Delta\Ss_R=c(\Ss_R\otimes\Ss_R)\Delta.$$
If $\Ss$ is invertible with respect to composition, then
$$\St_R m=m(\St_R\otimes\St_R)c^{-1}, \qquad
        \Delta\St_R=c^{-1}(\St_R\otimes\St_R)\Delta.$$
\end{lem}
\begin{pf}
Since $m$ is a coalgebra morphism and $\Ss$ is the inverse of the identity
in the monoid $\Hom_{\cate}(R,R)$, we have that $\Ss m$ is the inverse
of $m$ in the monoid $\Hom_{\cate}(R\otimes R,R)$. We have moreover that 
\begin{align*}
m*(m(\Ss\otimes\Ss)c)&= m (m\otimes (m(\Ss\otimes\Ss)c))
        (\id\otimes c\otimes\id)(\Delta\otimes\Delta) \\
&=m(m\otimes m)(\id\otimes\id\otimes\Ss\otimes\Ss)(\id\otimes c_{R,R\otimes R})
        (\Delta\otimes\Delta) \\
&=m(\id\otimes m)(\id\otimes m\otimes\id)(\id\otimes\id\otimes\Ss\otimes\Ss)
        (\id\otimes\Delta\otimes\id)(\id\otimes c)(\Delta\otimes\id) \\
&=m(m\otimes\id)(\id\otimes\Ss\otimes\id)(\Delta\otimes u\cou) \\
&=u\cou\otimes u\cou=u_{R\otimes R}\cou_{R\otimes R}.
\end{align*}
whence the first equality.

The second equality in the first line follows dualizing the equality just
proved. The second line follows immediately from the first using the naturality
of $c$.
\end{pf}

\bigskip

\begin{rem}
Let us suppose that there exists a forgetful functor
$U:\cate\to\cate'$ into some monoidal category
$\cate'$ in such a way that if $U(f)$ is an isomorphism
then $f$ is an isomorphism. Let $H$ be a bialgebra in $\cate$.
$UH$ is then an algebra and a coalgebra in $\cate'$ (in general
it is not a bialgebra). If there exists an antipode for $UH$ in
$\cate'$ (i.e., if the identity morphism has an inverse in the monoid
$\Hom_{\cate'}(UH,UH)$) then there exists an antipode for $H$ in $\cate$,
namely, $\Ss_{\cate}$ is the morphism such that 
$U(\Ss_{\cate})=\Ss_{\cate'}$, which exists by our hypothesis on $U$, as
we now prove. Consider the morphism in $\cate$
given by $F:H\otimes H\to H\otimes H,\ F=(\id\otimes m)(\Delta\otimes\id)$.
$U(F)\in\End_{\cate'}(UH\otimes UH)$ is an isomorphism in $\cate'$,
whose inverse is $(\id\otimes m)(\id\otimes\Ss_{\cate'}\otimes\id)
(\Delta\otimes\id)$. Let $T\in\End_{\cate}(H\otimes H)$
be the inverse of $F$. The antipode is then given by the composition
$$\Ss_{\cate}=\left(H @>r_H^{-1}>> H\otimes\uni @>\id\otimes u>> H\otimes H
        @>T>> H\otimes H @>\cou\otimes\id>> \uni\otimes H @>l_H>> H\right).$$
In the usual cases, $\cate'$ is the category of $\k$-vector spaces.
The hypothesis is verified for instance when $\cate=\yd$ for $H$ a
Hopf algebra with bijective antipode, or $\cate=\modui H$
for $H$ a QT bialgebra.
\end{rem}

Let us recall the definition of $\Im M$ for $M\in\yd$ stated after the
definition \ref{defyd}. It is immediate to see with the previous remark
(or by direct computation) that if $R$ is a Hopf algebra in $\yd$ then
$\Im R$ is also a Hopf algebra.

\subsection{Bosonization}\ \\
\label{bosoni}
Let now $H$ be a fixed Hopf algebra over $\k$ with bijective antipode.
There is a one-to-one correspondence between Hopf algebras in $\yd$ and Hopf
algebras $A$ with morphisms of Hopf algebras
$$A \mathop{\leftrightarrows}^{\iota}_p H$$
such that $p{\iota}=\id_H$. This correspondence was found by Radford in
\cite{radford:85} and explained in these terms by Majid in \cite{majid:94a}.
We give the details here:

Let $\displaystyle A \mathop{\leftrightarrows}^{\iota}_p H$ be as above. Let
$R=\coind AH=\LKer p=\{a\in A\ |\ (\id\otimes p)\Delta(a)=a\otimes 1\}$.
It is immediate
that this is a subalgebra of $A$, with the same unit. The counit of $R$ is
the restriction of that of $A$. We define the comultiplication, the antipode,
the action and the coaction by
\begin{align*}
\Delta_R(r)&=\com r1 ({\iota}\Ss_H(p\com r2))\otimes\com r3, \\
\Ss_R(r)&=({\iota}p(\com r1))\Ss_A(\com r2), \\
h\ganda r&=\com h1r\Ss\com h2, \\ 
\delta(r)&=(p\otimes\id)\Delta(r).
\end{align*}
It is straightforward to see that these morphisms make $R$ into a Hopf
algebra in $\yd$.

\medskip

Conversely, if $R$ is a Hopf algebra in $\yd$, let $A=R\# H$ be the
semidirect product algebra build from the action of $H$ on $R$, and let
\begin{align*}
&\Delta_A(r\# h)=(\com r1\#\comn r21\com h1)\otimes(\comm r20\#\com h2), \\
&{\iota}(h)=1\# h,\qquad p(r\# h)=\cou(r)h, \\
&\Ss_A(r\# 1)={\iota}(\Ss(\con r1))(\Ss_R(\com r0)\# 1)
	=(\Ss(\con r1)\ganda\Ss_R(\com r0)\# \Ss(\con r2)).
\end{align*}
These morphisms make $R\# H$ into a Hopf
algebra, and the constructions are mutually inverse.
Majid calls $R\# H$ the ``bosonization" of $R$.

\subsection{Examples of braided Hopf algebras}
\label{ejem}
\addtocounter{subsubsection}{-1}
\subsubsection{}
Let $H=\k$. The Yetter--Drinfeld category over $H$ reduces
in this case to the category of vector spaces over $\k$ (with trivial
actions and coactions), and the braiding is just the usual flip
$x\otimes y\mapsto y\otimes x$. A Hopf algebra in this category is just a
(classic) Hopf algebra over $\k$.

\subsubsection{}
Let $N$ be a natural number, $\xi$ a primitive $N$-root of unity in $\k$
and $A=T_{\xi,N}$ the Taft algebra of order $N^2$ over $\k$, which is
generated as a $\k$ vector space by the elements
$\{g^ix^j\}_{0\le i,j \le N-1}$, with relations $g^N=1,\ x^N=0$, 
and $xg=\xi gx$. The comultiplication is given by $\Delta g=g\otimes g$,
and $\Delta x=g\otimes x+x\otimes 1$. The antipode is given by $\Ss g=g^{-1}$,
and $\Ss x=-g^{-1}x$. The counit, by $\cou g=1,\ \cou x=0$.

Let $H$ be the group algebra of the cyclic group of $N$ elements. We shall
denote also by $g$ a generator of this group. There is a morphism of Hopf
algebras
$$\pi:A\to H,\quad \pi(g^ix^j)=g^i\delta_{j,0}.$$
This morphism has a section, namely
$$\iota:H\to A,\quad \iota(g^i)=g^i.$$
One sees that $R=\coind AH$ is isomorphic as an algebra to
$\k[x]/(x^N)$. It has comultiplication $\Delta_R x=x\otimes 1+1\otimes x$,
counit $\cou_R(x)=0$, and antipode $\Ss_Rx=-x$. 
The action and coaction of $H$ over $R$ are given by $g\ganda x=gxg^{-1}
=\xi^{-1}x$, and $\delta(x)=g\otimes x$.

\subsubsection{}
Let $H$ be as before the group algebra of a cyclic group of order $N$ with
generator $g$.\linebreak
Let $A={\h(\xi,m)}={\k\gp{y,x,g}/\sim}$ be the book algebra considered
in \cite{andrus-schneider:98}.
It is the Hopf algebra with generators $\{x,y,g\}$ and relations
$$x^N=y^N=0,\quad g^N=1,\quad,gx=\xi xg,\quad gy=\xi^{m}yg,\quad xy=yx$$
and with comultiplication, antipode and counit given by
\begin{gather*}
\Delta(x)=x\otimes g+1\otimes x,\quad \Delta(y)=y\otimes 1+g^m\otimes y,\quad
\Delta(g)=g\otimes g,\\
\Ss(x)=-xg^{-1},\quad \Ss(y)=-g^{-m}y,\quad
\Ss(g)=g^{-1},\quad \cou(x)=\cou(y)=0,\quad \cou(g)=1.
\end{gather*}

One can take here either ${H=T_{\xi,N}}={\k\gp{x,g}/\sim}$
or $H=\k\gp{g}$. In the first case, $p(y)=0,\ p(x)=x,\ p(g)=g$, and
$$\LKer(p)=\frac{\k[y]}{(y^N)}.$$
In the second case, let $\bar x=xg^{-1}$, and $p(y)=p(x)=0,\ p(g)=g$.
Then $$\LKer(p)=\frac{\k\gp{\bar x,y}}{(\bar{x}^N,y^N,\bar xy-\xi^my\bar x)}.$$

\subsubsection{}
The preceding examples are particular cases of a wider class of
braided Hopf algebras, which we now define.
Suppose $\Gamma$ is an abelian group.
Let $g_1,\ldots,g_n\in\Gamma$ and $\chi_1,\ldots,\chi_n:\Gamma\to\k^{\times}$
characters. Suppose that for $i\neq j$ we have $\chi_i(g_j)\chi_j(g_i)=1$.
Let $N_i$ be the order of $\chi_i(g_i)$, and $q_{i,j}=\chi_j(g_i)$.
Let $R$ be the algebra generated
by elements $x_1,\ldots,x_n$ with relations
\begin{alignat*}{2}
x_i^{N_i}&=0\quad &&\forall i, \\
x_ix_j&=q_{i,j}x_jx_i \quad &&\hbox{ if } i\neq j.
\end{alignat*}
Thus the set of monomials $\{x_1^{r_1}\cdots x_n^{r_n},\ 0\le r_i\le N_i\}$
is clearly a basis for $R$. We define the action and coaction of $H$ by
\begin{align*}
g\ganda(x_1^{r_1}\dots x_n^{r_n}) &= \chi_1(g)^{r_1}\dots\chi_n(g)^{r_n}
                                        x_1^{r_1}\dots x_n^{r_n}, \\
\delta(x_1^{r_1}\dots x_n^{r_n}) &= g_1^{r_1}\dots g_n^{r_n}\otimes
                                        x_1^{r_1}\dots x_n^{r_n},
\end{align*}
and the comultiplication, counit and antipode by
\begin{align*}
\Delta(x_i) &=1\otimes x_i+x_i\otimes 1, \\
\cou x_i &=0, \\
\Ss x_i &=-x_i.
\end{align*}
Then $R$ is a braided Hopf algebra over $\k\Gamma$. Following Manin,
these Hopf algebras are called {\em quantum linear spaces}.
Several classification results were obtained in \cite{andrus-schneider:97} from
the study of these braided Hopf algebras, including the classification
of pointed Hopf algebras of order $p^3$, $p$ an odd prime.

\subsubsection{}
Let $\g$ be a complex finite dimensional simple Lie algebra, $\bo$ a
Borel subalgebra, $A$ the Cartan matrix of $\g$. The Lusztig's algebras
$\f$ and ${}'\f$ constructed from $A$ are braided Hopf algebras in
$\yd$, where $H$ is the group algebra of a free abelian group.
See example \ref{ejelus}, or \cite{schauenburg:96}, \cite{lusztig:93}
for the details. The bosonization of $\f$ is the
quantized enveloping algebra $\U_q(\bo)$ of $\bo$. Since Drinfeld
showed how to obtain the quantized enveloping algebra $\U_q(\g)$
of $\g$ from $\U_q(\bo)$ via the double construction, we see that
quantum groups can be derived in a conceptual way from the setting
of braided Hopf algebras.

\section{Duals and opposite algebras. Integrals}
\subsection{First results about duals}
\begin{prop} \label{ydesri}
Let $H$ be a Hopf algebra over $\k$ with bijective antipode.
Then $\yd$ is a braided rigid category.
\end{prop}

\begin{pf}
We have seen that $\yd$ is braided. We shall prove now rigidity.
Let $M$ be an object of $\yd$. Let $\{\bdm\}_{\alpha\in A}$ be
a basis of $M$ as a $\k$-vector space, and $\{\dbm\}_{\alpha\in A}$
its dual basis. We shall omit the summation symbol in any formula
with the occurrence of the element $\sum_{\alpha\in A}\bdm\otimes\dbm$.
We take $\dai M$ and $\dad M$ to be the dual of $M$ as $\k$-vector spaces,
with the following structure:
\begin{equation*}
\left.\begin{split}
(h\cdot f)(m)&=f(\Ss(h)m), \\
\con f1 \otimes\com f0&=\St(\con{\bdm}1)\otimes f(\com{\bdm}0)\dbm
\end{split}\right\}\quad\mbox{for }\dad M,
\end{equation*}
and
\begin{equation*}
\left.\begin{split}
(h\cdot f)(m)&=f(\St(h)m), \\
\con f1\otimes\com f0&=\Ss(\con{\bdm}1)\otimes f(\com{\bdm}0)\dbm
\end{split}\right\}\quad\mbox{for }\dai M.
\end{equation*}

The morphisms $br,\ bl,\ dr,\ dl$ are the canonical morphisms
\begin{alignat*}{2}
br=\iota&:\k\to M\otimes\dad M \qquad && 1\mapsto\bdm\otimes\dbm, \\
bl=\iota&:\k\to\dai M\otimes M \qquad && 1\mapsto\dbm\otimes\bdm, \\
dr=\ev&:\dad M\otimes M\to\k \qquad && f\otimes m\mapsto f(m), \\
dl=\ev&:M\otimes \dai M\to\k \qquad && m\otimes f\mapsto f(m),
\end{alignat*}
which are morphisms in $\yd$ with respect to the above defined structures.
It is immediate that these morphisms satisfy equations \eqref{eqdual}.
\end{pf}

\bigskip

Let $\cate$ be any braided rigid category. We recall that this means
(for us) that there exist functors $M\mapsto\dad M$ and $M\mapsto\dai M$.
These functors are inverse to each other via canonical isomorphisms
since $\dai{(\dad M)}\simeq \dad{(\dai M)}\simeq M$.
Indeed, it is clear that
$(M,\,br_M:\uni\to M\otimes\dad M,\,dr_M:\dad M\otimes M\to\uni)$
satisfies the axioms of left dual for $\dad M$, and
$(M,\,bl_M:\uni\to\dai M\otimes M,\,dl_M:M\otimes\dai M\to\uni)$
satisfies the axioms of right dual for $\dai M$.
Moreover, the functors $M\mapsto\dad M$ and $M\mapsto\dai M$ are naturally
isomorphic, as can be seen considering
\begin{gather*}
\dad M @>\id\otimes bl>> \dad M\otimes \dai M\otimes M @>c\otimes\id>>
	\dai M\otimes\dad M\otimes M @>\id\otimes dr>> \dai M, \\
\dai M @>br\otimes\id>> M\otimes\dad M\otimes\dai M @>\id\otimes c^{-1}>>
	M\otimes\dai M\otimes\dad M @>dl\otimes\id>> \dad M.
\end{gather*}
Thus, $M^{**}\simeq {}^*(M^*)\simeq M\simeq (\dai M)^*\simeq {}^{**}M$
via natural isomorphisms.
In the category $\yd$ these morphisms are given by
\begin{alignat*}{2}
&M\to M^{**},\qquad && m\mapsto\Ss(\con m1)\com m0, \\
&M\to {}^{**}\!M,\qquad && m\mapsto\Ss^{-2}(\con m1)\com m0.
\end{alignat*}
The rather strange asymmetry between both morphisms comes from the
fact that we use $c^{-1}$ in the first one and $c$ in the second one.

\begin{defn}
Let $\cate$ be any rigid category, and $M,N$ be objects of $\cate$.
Let $F:M\to N$ be a morphism in $\cate$. We define the transposes of $F$
as
\begin{align*}
\dad F &= \dad N \to \dad N\otimes\uni \to \dad N\otimes M\otimes\dad M
	@>\id\otimes F\otimes\id>> \dad N\otimes N\otimes\dad M \to
	\uni\otimes\dad M \to \dad M, \\
\dai F &= \dai N \to \uni\otimes\dai N \to \dai M\otimes M\otimes\dai N
	@>\id\otimes F\otimes\id>> \dai M\otimes N\otimes\dai N \to
	\dai M\otimes\uni \to \dai M. 
\end{align*}
\end{defn}

\begin{rem}
\label{F*}
Most of the rigid categories we consider are subcategories of the
category of $\k$-vector spaces, and the duals are preserved by the
forgetful functor (the maps $br,\ bl,\ dr,\ dl$ are also preserved).
When this happens, the maps $\dad F$ and $\dai F$ coincide (via the
forgetful functor) with the usual transpose map of $F$.
We observe that this means that for $F:M\to N$ a morphism in $\yd$, the
transpose as $\k$-vector spaces $\dad F:\dad N\to \dad M$ is a morphism
in $\yd$.
\end{rem}

\begin{lem}
\label{*davuel}
Let $\cate$ be any rigid category, and let $M,N\in\cate$.
There exist natural isomorphisms
\begin{align*}
\dad\phi_{M,N} &: \dad M\otimes\dad N \to \dad {(N\otimes M)}, \\
\dai\phi_{M,N} &: \dai M\otimes\dai N \to \dai {(N\otimes M)}.
\end{align*}
\end{lem}
\begin{pf}
To prove that $\dad M\otimes\dad N \simeq \dad {(N\otimes M)}$ it would be
sufficient to prove that $\dad M\otimes\dad N$ satisfy \eqref{eqdual} for
certain morphisms $br,\ bl,\ dr$ and $dl$, but in order to prove naturality
it is necessary to give the explicit definition of $\dad\phi$ and $\dai\phi$.
\begin{multline*}
\dad\phi=\dad M\otimes\dad N @>\id\otimes\id\otimes br_{N\otimes M}>>
(\dad M\otimes\dad N)\otimes(N\otimes M)\otimes\dad{(N\otimes M)}\to \\
\dad M\otimes(\dad N\otimes N)\otimes\dad M\otimes\dad{(N\otimes M)}
@>\id\otimes dr_N\otimes\id\otimes\id>>
	\dad M\otimes M\otimes\dad{(N\otimes M)}
	@>dr_M\otimes\id>> \dad{(N\otimes M)}.
\end{multline*}
Analogously for $\dai\phi$. The proof that $\dai\phi$ and $\dad\phi$ are
natural is straightforward but tedious and we omit it.
\end{pf}

\begin{lem}\label{cphi}
Let $N$, $M$ be objects of $\cate$. We have
$$c^*_{M,N}\dad\phi_{M,N}=\dad\phi_{N,M}c_{\dad M,\dad N}.$$
\end{lem}
\begin{pf}
First, we claim that
\begin{multline*}
\uni \to (M\otimes N)\otimes\dad{(M\otimes N)}
@>c\otimes\id>>(N\otimes M)\otimes\dad{(M\otimes N)} \\
=\uni \to (N\otimes M)\otimes\dad{(N\otimes M)}
@>\id\otimes\dad c>>(N\otimes M)\otimes\dad{(N\otimes M)}
\end{multline*}
In fact, tensoring both sides with $(M\otimes N)$ on the right and composing
with $dr_{(M\otimes N)}$ one gets $c$, whence the claim.

Second, we claim that
\begin{multline*}
\dad M\otimes\dad N\otimes M\otimes N @>\id\otimes c>>
	\dad M\otimes\dad N\otimes N\otimes M \to \dad M\otimes M\to\uni \\
=\dad M\otimes\dad N\otimes M\otimes N @>c\otimes\id>>
	\dad N\otimes\dad M\otimes M\otimes N\to\dad N\otimes N\to\uni.$$
\end{multline*}
In fact, both sides equal
$\dad M\otimes\dad N\otimes M\otimes N \flllad{\id\otimes c^{-1}\otimes\id}
	\dad M\otimes M\otimes\dad N\otimes N\to\uni.$
Thus,
\begin{multline*}
\dad\phi_{N,M}c_{\dad M,\dad N} \\
\shoveleft{=\dad M\otimes\dad N @>c_{\dad M,\dad N}>>
\dad N\otimes\dad M \to \dad N\otimes\dad M \otimes (M\otimes N)
	\otimes\dad{(M\otimes N)} @>dr\otimes\id>> \dad{(M\otimes N)}} \\
\shoveleft{=\dad M\otimes\dad N \to
\dad M\otimes\dad N \otimes (M\otimes N) \otimes\dad{(M\otimes N)}
@>c_{\dad M,\dad N}\otimes\id\otimes\id>>} \\
\shoveright{\dad N\otimes\dad M \otimes (M\otimes N) \otimes\dad{(M\otimes N)}
@>dr\otimes\id>>
\dad{(M\otimes N)}} \\
\shoveleft{=\dad M\otimes\dad N \to
\dad M\otimes\dad N \otimes (M\otimes N) \otimes\dad{(M\otimes N)}
@>\id^{\otimes 2}\otimes c_{M,N}\otimes\id>>} \\
\shoveright{\dad M\otimes\dad N \otimes (N\otimes M) \otimes\dad{(M\otimes N)}
@>dr\otimes\id>>
\dad{(M\otimes N)}} \\
\shoveleft{=\dad M\otimes\dad N \to
\dad M\otimes\dad N \otimes (N\otimes M) \otimes\dad{(N\otimes M)}
@>\id\otimes\id\otimes\dad{(c_{M,N})}>>} \\
\shoveright{\dad M\otimes\dad N \otimes (N\otimes M) \otimes\dad{(M\otimes N)}
@>dr\otimes\id>>
\dad{(M\otimes N)}} \\
\shoveleft{=\dad M\otimes\dad N \to
\dad M\otimes\dad N \otimes (N\otimes M) \otimes\dad{(N\otimes M)}
@>dr\otimes\id>>
\dad{(N\otimes M)} @>\dad{(c_{M,N})}>>\dad{(M\otimes N)} } \\
=\dad{(c_{M,N})}\dad\phi_{M,N}.
\end{multline*}
\end{pf}

\subsection{Equivalence of some Yetter--Drinfeld categories
	and dual Hopf algebras} \ \\
In the setting of Yetter-Drinfeld categories, one often needs to pass
from one category to another. This is usually possible. We give the
corresponding functors. We recall from \ref{rmfunmon} the definition
of a functor between braided categories.

\begin{prop}
\label{cateigua}
\begin{enumerate}\item Let $H$ be a Hopf algebra with bijective
antipode. The following categories are equivalent as braided
categories:
$$(i)\ \yd,\quad (ii)\ \hbyd,
	\quad (iii)\ \ydh,\quad (iv)\ \ydhb.$$

\item If $H$ is finite dimensional the preceding categories are
equivalent to the following ones (as braided categories):
$$(v)\ \ydhs,\quad (vi)\ \ydhsb,
	\quad (vii)\ \hsyd,\quad (viii)\ \hsbyd.$$

\item The following categories are equivalent as braided categories
if $H$ is finite dimensional:
$$(ix)\ \ydb,\quad (x)\ \ydhc.$$
\end{enumerate}
\end{prop}

\begin{pf}
For (1) and (2), let $M$ be an object in $\yd$. We first prove that
$(i)$, $(ii)$, $(v)$, $(vi)$ are equivalent. We take the structure
\begin{alignat*}{3}
& h\ganda^2m=\Ss(h)m,\quad && \delta^2(m)=\St\con m1\otimes\com m0\quad
	&&\mbox{for }\hbyd, \\
& m\gania^5f=\com m0\gi f,\con m1\gd,\quad && \delta^5(m)=
	\bdh m\otimes\dbh\quad &&\mbox{for }\ydhs, \\
& m\gania^6f=\com m0\gi f,\St\con m1\gd,\quad && \delta^6(m)=
	\bdh m\otimes\Ss(\dbh)\quad &&\mbox{for }\ydhsb.
\end{alignat*}
It is not difficult to verify that these structures make $M$ into
objects in the stated categories and that preserve tensor products.
In all these cases the natural isomorphism $\eta$ of Remark
\ref{rmfunmon} is the identity, i.e. $F(M\otimes N)=FM\otimes FN$.

We verify the compatibility with the braiding between $\yd$ and $\ydhs$.
The others are analogous.
Let $c_1$ and $c_5$ denote the braidings in the
respective categories. Let $M$ and $N$ be objects in $\yd$. Then we have
to prove that $Fc_1=c_5F:M\otimes N\to N\otimes M$. Let
$m\otimes n\in M\otimes N$, and denote by the same symbol the corresponding
element in $FM\otimes FN=F(M\otimes N)$. Then
\begin{align*}
c_5(m\otimes n) &=\delta^3_0n\otimes m\delta^3_1n=
	\bdh n\otimes m\gania^3\dbh \\
&=\bdh n\otimes\com m0\gi\dbh,\con m1\gd=
	\con m1 n\otimes\com m0=c_1(m\otimes n).
\end{align*}

In an analogous way it can be proved that the categories
$(iii)$, $(iv)$, $(vii)$ and $(viii)$ are equivalent. We give
the equivalence between $(i)$ and $(iii)$, which is more subtle
since $\eta\neq\id$.
Let $M$ be an object in $\yd$. We define $\Re(M)$ in $\ydh$ to be
$M$ as a $\k$-vector space, with the structure given by
$$m\gania^3h=\St(h)m,\quad \delta^3m=\com m0\otimes\Ss\con m1,$$
whence
$$\delta^3(m\gania^3h)=
\delta^3_0m\gania^3\com h2\otimes\Ss(\com h1)\delta^3_1(m)\com h3.$$
Observe that there is a natural isomorphism
$$\phi=\phi_{M,N}:\Re(M\otimes N)\to \Re N\otimes\Re M,\quad
	(m\otimes n)\mapsto n\otimes m.$$
We define
$$\eta_{M,N}=\phi_{M,N}^{-1}\circ c_{\Re M,\Re N}
	:\Re(M)\otimes\Re(N)\to\Re(M\otimes N),$$
that is,
\begin{align*}
\eta_{M,N}(m\otimes n)
	&=\phi_{M,N}^{-1}(\delta^5_0(n)\otimes m\gania^5\delta^5_1(n)) \\
	&=\phi_{M,N}^{-1}(\com n0\otimes\St(\Ss\com n1)m) \\
	&=\phi_{M,N}^{-1}(\com n0\otimes\com n1m)=\com n1m\otimes\com n0
		\in\Re(M\otimes N).
\end{align*}
It is straightforward to check that $(\Re,\eta)$ is a functor between braided
categories. We verify for instance \eqref{funmon}:
\begin{align*}
\eta_{M\otimes N,P}\circ(\eta_{M,N}\otimes\id)(m\otimes n\otimes p)
	&=\eta_{M\otimes N,P}(\con n1m\otimes\com m0\otimes p) \\
	&=\con p1(\con n1m\otimes\com m0)\otimes\com p0 \\
	&=\con p2\con n1m\otimes\con p1\com n0\otimes\com p0, \\
\eta_{M,N\otimes P}\circ(\id\otimes\eta_{N,P})(m\otimes n\otimes p)
	&=\eta_{M,N\otimes P}(m\otimes\con p1n\com p0) \\
	&=\con p4\con n1\Ss(\con p2)\con p1m\otimes\com p3\com n0
		\otimes\com p0 \\
	&=\con p2\con n1m\otimes\con p1\com n0\otimes\com p0.
\end{align*}

(3) Let $M$ be an object in $\ydb$, and define
$$m\gania^{10}h=\St(h)m,\quad \delta^{10}(m)=\com m0\otimes\con m1.$$
As before, it is straightforward to see that this is an object in $\ydhc$,
and that the braiding is that of $\ydb$.
\end{pf}

\bigskip

We concentrate now on dual Hopf algebras. It would be possible to define
the dual of a braided Hopf algebra in $\yd$ declaring $\dad R$
(resp. $\dai R$) to be the right dual (resp. left dual) of $R$ in
$\yd$ with the algebra and coalgebra structure transposes of the
coalgebra and algebra structures of $R$. This would fail to be a
bialgebra in $\yd$ because the compatibility between multiplication
and comultiplicatin does not transpose to the compatibility between
the transpose operations. There are two ways to fix this problem.
The first one is to take a kind of $\rsbop$, which is a Hopf algebra
in $\yd$ (as it is done for a general rigid braided category in
\cite{majid:94b}). The second one is to consider $\dad R$ (or $\dai R$)
as a Hopf algebra in $\hsyd$ and recover a Hopf algebra in $\yd$ via
$\ba\Re$, the inverse functor to $\Re$.
The natural way to see $\dad R$ as a Hopf algebra in $\hsyd$
is by means of the following construction: let
$$\displaystyle R\# H \mathop{\leftrightarrows}^{\iota}_p H$$
be the construction given in section \ref{bosoni}. We can dualize it to get
$$R^*\otimes H^*\simeq (R\# H)^*
	\mathop{\rightleftarrows}^{{\iota}^*}_{p^*} H^*,$$
where the first is an isomorphism of vector spaces, and
${\iota}^*p^*=\id_{H^*}$. It is immediate that
$\LKer({\iota}^*)=\dad R\otimes\cou_H\subseteq\dad R\otimes\dad H$.
This makes $R^*$ into a Hopf algebra in $\hsyd$. One can get $\dai R$
with the same procedure, but starting out with $\Im R$ instead of $R$.

\medskip

We prefer instead of doing this the more categorical way: we define
duals of a Hopf algebra in any rigid braided category as it is done
by several authors (see for instance \cite{takeuchi:97}).
For the special case of $\yd$, we get the above duals.

\begin{defn} \label{defdualcat}
Let $\cate$ be any braided rigid category and $M,N$ objects of
$\cate$. Let $\dad\phi$ and $\dai\phi$ be the isomorphisms of
\ref{*davuel}. We define
$$\dad\sigma_{M,N} = \dad M\otimes\dad N \fllad{c^{-1}}
	\dad N\otimes\dad M \fllad{\dad\phi} \dad{(M\otimes N)},$$
and then we define the structure of $\dad R$ by
\begin{align*}
m_{\dad R}&=\dad R\otimes\dad R @>\dad\sigma>> \dad {(R\otimes R)}
	@>\dad\Delta>> \dad R, \\
\Delta_{\dad R}&=\dad R @>\dad m>> \dad {(R\otimes R)}
	@>(\dad\sigma)^{-1}>> \dad R\otimes \dad R, \\
\Ss_{\dad R}&= \dad{(\Ss_R)},
\ u_{\dad R}=\dad{(\cou_R)},\ \cou_{\dad R}=\dad{(u_R)}.
\end{align*}
We define $\dai R$ in the same manner, replacing the duals on the
right by duals on the left.
\end{defn}

\begin{lem}
These morphisms make $\dad R$ and $\dai R$ into Hopf algebras in $\cate$.
\end{lem}
\begin{pf}
We use lemma \ref{cphi}. It is straightforward to prove associativity,
coassociativity and the axioms for unit, counit and antipode. We shall
prove the compatibility between multiplication and comultiplication
for $\dad R$. The proof for $\dai R$ is analogous. Let $M,\ N,\ S$ and
$T$ be objects in $\cate$. We denote 
\begin{align*}
c_{2,2}&:(M\otimes N\otimes S\otimes T)\to (S\otimes T\otimes M\otimes N)
	=c_{M\otimes N,S\otimes T}, \\
\dad\phi_{2,2}&:\dad{(M\otimes N)}\otimes\dad{(S\otimes T)}\to
	\dad{(S\otimes T\otimes M\otimes N)}
	=\dad\phi_{M\otimes N,S\otimes T}, \\
\dad\phi_4&:(\dad M\otimes\dad N\otimes\dad S\otimes\dad T)\to
	\dad{(T\otimes S\otimes N\otimes M)}
	=\dad\phi_{2,2}(\dad\phi\otimes\dad\phi).
\end{align*}
Let us observe that if $f$ and $g$ are morphisms then 
$$(\dad f\otimes\dad g)={(\dad\phi)}^{-1}\dad{(g\otimes f)}\dad\phi.$$
We claim that
$$c(\dad\phi_{2,2})^{-1}\dad{(\id\otimes c\otimes\id)}\dad\phi_{2,2}c^{-1}
	=(\dad\phi\otimes\dad\phi)(c^{-1}\otimes c^{-1})
	(\id\otimes c\otimes\id)(c\otimes c)
	((\dad\phi)^{-1}\otimes(\dad\phi)^{-1}).$$
This is true because
\begin{align*}
\dad\phi_{2,2} &c^{-1}(\dad\phi\otimes\dad\phi)(c^{-1}\otimes c^{-1})
	(\id\otimes c\otimes\id)(c\otimes c)
	({(\dad\phi)}^{-1}\otimes{(\dad\phi)}^{-1})c{(\dad\phi_{2,2})}^{-1} \\
&=\dad\phi_{2,2}(\dad\phi\otimes\dad\phi)[c^{-1}_{2,2}(c^{-1}\otimes c^{-1})
	(\id\otimes c\otimes\id)(c\otimes c)c_{2,2}]
	({(\dad\phi)}^{-1}\otimes{(\dad\phi)}^{-1}){(\dad\phi_{2,2})}^{-1} \\
&=\dad\phi_4(\id\otimes c\otimes\id){(\dad\phi_4)}^{-1}=
	\dad{(\id\otimes c\otimes\id)}.
\end{align*}
Hence
\begin{align*}
\Delta_{\dad R}m_{\dad R}
&= c(\dad\phi)^{-1}\dad m\dad\Delta\dad\phi c^{-1} 
	=c(\dad\phi)^{-1}\dad{(\Delta m)}\dad\phi c^{-1} \\
&=c{(\dad\phi)}^{-1}\dad{\left[(m\otimes m)(\id\otimes c\otimes\id)
	(\Delta\otimes\Delta)\right]}\dad\phi c^{-1} \\
&=c{(\dad\phi)}^{-1}\dad{(\Delta\otimes\Delta)}\dad{(\id\otimes c\otimes\id)}
	\dad{(m\otimes m)}\dad\phi c^{-1} \\
&=c(\dad\Delta\otimes\dad\Delta){(\dad\phi_{2,2})}^{-1}\dad{(\id\otimes c\otimes\id)}
	\dad\phi_{2,2}(\dad m\otimes\dad m)c^{-1} \\
&=(\dad\Delta\otimes\dad\Delta)c{(\dad\phi_{2,2})}^{-1}\dad{(\id\otimes c\otimes\id)}
	\dad\phi_{2,2}c^{-1}(\dad m\otimes\dad m) \\
&=(\dad\Delta\otimes\dad\Delta)(\dad\phi\otimes\dad\phi)(c^{-1}\otimes c^{-1})
	(\id\otimes c\otimes\id)(c\otimes c)
	({(\dad\phi)}^{-1}\otimes{(\dad\phi)}^{-1})(\dad m\otimes\dad m) \\
&=(m_{\dad R}\otimes m_{\dad R})(\id\otimes c\otimes\id)
	(\Delta_{\dad R}\otimes\Delta_{\dad R})
\end{align*}
\end{pf}

Let $R$ be a Hopf algebra in $\yd$. We give the specific structure
for $\dad R$ and $\dai R$. The formulae are exactly the same for
both algebras.
\begin{align*}
\gi m(f\otimes g),r\gd &=
	\gi f,\comm r20\gd\gi g,\St(\comn r21)\com r1 \gd \\
	&= \gi\com g0,\com r2\gd\gi\St(\con g1)f,\com r1\gd. \\
\gi \com f1,r\gd\gi\com f2,s\gd &=
	\gi f,mc(s\otimes r)\gd = \gi f,(\con s1r)\com s0\gd \\
	&= \gi\St(\comn f21)\com f1,\con s1r\gd\gi\comm f20,\com s0\gd. \\
\gi\Ss f,r\gd &= \gi f,\Ss r\gd,\quad
	\gi 1_{\dad R},r\gd=\gi\cou_R,r\gd,\quad
	\gi \cou_{\dad R},f\gd=\gi f,1_R\gd.
\end{align*}

The following result was found by many authors, see for instance
\cite{takeuchi:97} or \cite{bespa-ker-lyuba-tura}.
\begin{prop}
Let $\cate$ be a braided monoidal category. As usual, we denote by
$\cateb$ the same category but with the inverse braiding, i.e.
$c^{\cateb}_{M,N}=(c^{\cate}_{N,M})^{-1}$.
Let $R$ be a Hopf algebra in $\cate$ whose antipode is an isomorphism.
We define $\rop,\ \rcop$ and $\rbop$ by
\begin{alignat*}{3}
&m_{\rop}=m_R\circ c^{-1}_{R,R}, \qquad && \Delta_{\rop}=\Delta_R,
	&& \Ss_{\rop}=\St_R, \\
&m_{\rcop}=m_R, && \Delta_{\rcop}=c^{-1}_{R,R}\circ\Delta_R,\qquad
	&& \Ss_{\rcop}=\St_R, \\
&m_{\rbop}=m_R\circ c_{R,R}, && \Delta_{\rbop}=c^{-1}_{R,R}\circ\Delta_R,
	&& \Ss_{\rbop}=\Ss_R.
\end{alignat*}
and the other structure morphisms remain equal as those of $R$.
Then $\rop$ and $\rcop$ are Hopf algebras in $\cateb$, and $\rbop$
is a Hopf algebra in $\cate$.
\end{prop}

\begin{pf}
The general proof is straightforward (and in fact very easy using
drawings). We give a direct proof for the particular case of a
Yetter--Drinfeld category. We prove the statement for $\rop$.
The proof for $\rcop$ is analogous, and for $\rbop$ is the
composition of the other two. Associativity is easy to prove. We check
the compatibility between the multiplication and comultiplication:
\begin{align*}
\Delta \mulop(r\otimes s)
&=\Delta(\com s0(\St(\con s1)r)) \\
&=\comm s01\left(\commn s021\St(\con s1)\com r1\right)\otimes
	\commm s020\left(\St(\con s2)\com r2\right) \\
&=\comm s10\left(\comn s21\St(\comn s11\comn s22)\com r1\right)\otimes
	\comm s20\left(\St(\comn s12\comn s23)\com r2\right) \\
&=\comm s10\left(\St(\comn s11)\com r1\right)\otimes
	\comm s20\left(\St(\comn s21)\St(\comn s12)\com r2\right) \\
&=(\mulop\otimes \mulop)\left(\com r1\otimes
	\comm s10\otimes\St(\comn s11)\com r2\otimes\com s2\right) \\
&=(\mulop\otimes \mulop)(\id\otimes c^{-1}\otimes\id)
	(\Delta\otimes\Delta)(r\otimes s).
\end{align*}

It is straightforward, using \ref{sdavuel}, to check that $\St$ verifies
the axioms for the antipode.
\end{pf}

\begin{rem}
We note that the above definitions can be made for algebras or coalgebras
in a braided category. Then, if $A$ is an algebra in $\cate$ it can be
defined $\aop$ as the same object as $A$ with multiplication
$m_{\aop}=m_Ac^{-1}_{A,A}$, and if $C$ is a coalgebra in $\cate$ it can
be defined $\ccop$ as the same object as $C$ with comultiplication
$\Delta_{\ccop}=c^{-1}_{C,C}\circ\Delta_C$.
\end{rem}

\subsection{Integrals}\ \\
By classical results a finite dimensional Hopf algebra has a one dimensional
space of left (resp. right) integrals. These results can be generalized
to finite Hopf algebras in braided categories, as it is done in
\cite{doi:unp}, \cite{lyubashenko:95a,lyubashenko:95b} or \cite{takeuchi:97}
(we define an object in $\cate$ to be {\em finite} if it has a (right) dual 
\footnote{Takeuchi defines {\em dual} in a weaker form: he calls $\dad M$
{\em the dual of $M$} if there exists a morphism $\dad M\otimes M\to\uni$
with the universal property that for every morphism $f:X\otimes M\to\uni$,
there exists a unique morphism $F:X\to\dad M$ such that
$f=\left(X\otimes M\fllad{F\otimes\id}\dad M\otimes M \flad{\ev}\uni\right)$.
He define an object to be {\em finite} if it has a dual in the sense of
\ref{duales}. It is easy to see that if $M$ has a dual $\dad M$ in the
sense of \ref{duales} then $\dad M$ is the dual in the sense of Takeuchi.
The terminology is consistent with the usual cases: for instance,
if $H$ is a QT bialgebra, every module in $\Modui H$ has a dual in the
sense of Takeuchi, but it is finite if and only if it is finite
dimensional.} in the sense of \ref{duales}).
In what follows $\cate$ shall be a braided category which has equalizers.
As we noted above, monoidal categories are the natural context to define
algebras and coalgebras. Given an algebra $A$ in a monoidal category it
is routine to define a (left) $A$ module: it is a pair $(M,\ganda)$ with
$M$ an object in the category and $\ganda:A\otimes M\to M$ which verifies
\begin{alignat*}{2}
&\ganda\circ(\id\otimes\ganda)=\ganda\circ(m\otimes\id):
A\otimes A\otimes M\to M\quad && \mbox{(associativity),} \\
&\ganda\circ(u\otimes\id)\circ l^{-1}_M=\id:M\to M && \mbox{(unitary).}
\end{alignat*}

\medskip

Analogously is defined a right $A$-module. If $C$ is a coalgebra in the
category we define in a dual fashion a (left) $C$-comodule to be a pair
$(M,\delta)$ with $M$ an object in the category and
$\delta:M\to C\otimes M$ which verifies
\begin{alignat*}{2}
& (\Delta\otimes\id)\circ\delta=(\id\otimes\delta)\circ\delta:
	M\to C\otimes C\otimes M\quad && \mbox{(coassociativity),} \\
& l_M\circ(\cou\otimes\id)\circ\delta=\id:M\to M && \mbox{(counitary).}
\end{alignat*}

\medskip

Analogously is defined a right $C$-comodule. It is routine also to define
a Hopf module in a braided category:
\begin{defn}
Let $R$ be a Hopf algebra in $\cate$.
A {\em left $R$-Hopf module} is a triple $(M,\ganda,\delta)$, where 
$(M,\ganda)$ is a left $R$-module, $(M,\delta)$ is a left $R$-comodule and
$\delta$ is a morphism of modules, where the structure of left $R$-module
of $R\otimes M$ is given by
$$R\otimes R\otimes M @>\Delta_R\otimes\id\otimes\id>>
	R\otimes R\otimes R\otimes M @>\id\otimes c\otimes\id>>
	R\otimes R\otimes R\otimes M @>m_R\otimes\ganda>> R\otimes M.$$
As in the classic case, $\delta$ is a module morphism iff $\ganda$ is a
comodule morphism, where the comodule structure of $R\otimes M$ is given by
$$R\otimes M @>\Delta_R\otimes\delta>>
	R\otimes R\otimes R\otimes M @>\id\otimes c\otimes\id>>
	R\otimes R\otimes R\otimes M @>m_R\otimes\id\otimes\id>>
	R\otimes R\otimes M.$$
The definition of right $R$-Hopf module is analogous.
\end{defn}

\begin{defn}
Let $C$ be a coalgebra in $\cate$ with a unit $u:\uni\to C$ which is a
coalgebra morphism, and $(M,\delta)$ be a left $C$-comodule.
We define the space of coinvariants by means of the equalizer
$$\coini MC=\equ(M\mathop{\rightrightarrows}^{\delta}_{u\otimes\id}
		C\otimes M).$$
\end{defn}
The fundamental theorem for Hopf modules can be modified to braided
categories. Specifically,

\begin{prop}
Let $\cate$ be a braided category, let $R$ be a Hopf algebra in $\cate$,
and let $\hopmoi R$ be the category of left $R$-Hopf modules.
Then $\hopmoi R$ is equivalent to $\cate$ via
\begin{align*}
V\in\cate &\rightarrow
	(R\otimes V,m_r\otimes\id,\Delta_R\otimes\id)\in\hopmoi R, \\
\coini MR\in\cate &\leftarrow M\in\hopmoi R.
\end{align*}
\end{prop}
\begin{pf}
Mimic \cite[Th 4.1.1]{sweedler:69}. See also \cite[3.3]{bespa-draba:95}
for the case when $\cate$ is a braided category with split idempotents.
\end{pf}

This is the first step to prove the existence of non zero (left) integrals
in a finite dimensional Hopf algebra, and is used by Takeuchi in
\cite{takeuchi:97} in the same way. We follow now his work. If $R$
is a finite Hopf algebra in $\cate$, we define the structure of right
$R$-Hopf module on $\dad R$ given by
\begin{align*}
\ganid &:=\dad R\otimes R \fllllad{\id\otimes\Ss\otimes br}
	\dad R\otimes R\otimes R\otimes\dad R \fllllad{\id\otimes c\otimes\id}
	\dad R\otimes R\otimes R\otimes\dad R \fllllad{\id\otimes m\otimes\id}
	\dad R\otimes R\otimes \dad R @>dr\otimes\id>> \dad R \\
\delta &:=\dad R \fllad{\id\otimes br} \dad R\otimes R\otimes\dad R
	\fllllad{\id\otimes\Delta\otimes\id}
	\dad R\otimes R\otimes R\otimes\dad R
	\flllllad{\id\otimes c^{-1}\otimes\id}
	\dad R\otimes R\otimes R\otimes\dad R
	\fllad{dr\otimes c} \dad R\otimes R.
\end{align*}
The proof that $(\dad R,\ganid,\delta)$ is an $R$-Hopf module is
straightforward.

From the other hand, let $A$ be a coalgebra in $\cate$,
which is also a finite object. Then $\dad A$ is an algebra in $\cate$,
with multiplication $\dad m$ as in \ref{defdualcat}.
Moreover, if $u:\uni\to A$ is a coalgebra map, then $(\dad A,\dad u)$
is an augmented algebra in $\cate$. If $(M,\ganda)$ is a left
$\dad A$-module in $\cate$, we define on $M$ a structure
of right $A$-comodule as follows:
$$\rho=M\flllad{br\otimes\id} A\otimes\dad A\otimes M
	\flllad{\id\otimes\ganda} A\otimes M \flad{c} M\otimes A.$$
Then the invariants of $M$ are defined by the equalizer
$\invi M{\dad A}=\coind MA=\equ(\rho,\id\otimes u)$.

Let now $\dad R$ act on $\dad R$ on the left by multiplication.
We define the integrals by
$$\intgl{\dad R}=\invi{\dad R}{\dad R}=\coind{(\dad R)}R.$$
Then the right coaction $\rho$ coincides with $\delta$, as can be seen tensoring
both morphism on the left with $\dad R$, and then composing with
$(dr\otimes\id)(\id\otimes c^{-1})$. The fundamental theorem on Hopf modules
gives then
$$\dad R\simeq\intgl{\dad R}\otimes R.$$
Changing $(R,\dad R)$ by $(\dai R,R)$ we get
$$R\simeq\intgl R\otimes\dai R,$$
and thus $\dad R\simeq\intgl{\dad R}\otimes\intgl R\otimes\dai R\simeq
\intgl{\dad R}\otimes\intgl R\otimes\dad R$.
If the category has coequalizers, we can apply to this isomorphism
the functor $-\otimes_{\dad R}\uni$ and we get 
$$\uni\simeq\intgl R\otimes\intgl{\dad R},$$
which means that $\intgl R$ is an invertible object in $\cate$.
It is clear that the above construction can be made analogously to get
right integrals
$$\intgr{\dai R}=\invd{\dai R}{\dai R}=\coini{(\dai R)}R,$$
such that $\dai R\simeq R\otimes\intgr{\dai R}$, $R\simeq\dad R\otimes\intgr R$.

\bigskip

From the invertibility of the space of integrals, it is possible to deduce the
existence of a distinguished grouplike in $\dad R$, which reflects the
action of $R$ on the right over $\intgl R$. See \cite{takeuchi:97}.

\bigskip

\medskip
We want to compute now the defining equation for the space of integrals
$\intgl R$ for $R$ a Hopf algebra in $\yd$. Let $\{\bdr\}$, $\{\dbr\}$
be dual bases for $R$ and $\dai R$. We have
\begin{align*}
\rho(x) &= c(\id\otimes m_R)(br_{\dai R}\otimes\id)(x)=
	c(\id\otimes m_R)(\dbr\otimes\bdr\otimes x) \\
	&= c(\dbr\otimes\bdr x)=\con{\dbr}1(\bdr x)\otimes\com{\bdr}0.
\end{align*}
We then have for $y\in R$
\begin{align*}
(\id\otimes y)(\rho x) &= \con{\dbr}1(\bdr x)\gi\com{\dbr}0,y\gd
	=\Ss(\con{\bdrb}1)(\bdr x)\gi\dbr,\com{\bdrb}0\gd\gi\dbrb,y\gd \\
	&=\Ss(\con y1)(\bdr x)\gi\dbr,\com y0\gd
	=\Ss(\con y1)(\com y0 x).
\end{align*}
Therefore $x\in\intgl R$ iff
\begin{equation}\label{rin1}
\Ss(\con y1)(\com y0 x)=(\id\otimes y)(\rho x)
	=(\id\otimes y)(x\otimes\cou)=\cou(y)x\quad\forall y\in R.
\end{equation}
Hence, if $x\in\intgl R$ we have
\begin{equation}\label{rin2}
yx=\con y2\Ss(\con y1)(\com y0 x)=\con y1\cou(\com y0)x
	=\cou(y)x\quad\forall y\in R.
\end{equation}
Conversely, it is immediate to see \eqref{rin2} implies \eqref{rin1}.
Thus, for $R$ a Hopf algebra in $\yd$ we have the well known equation
\begin{equation}\label{rinte}
x\in\intgl R \sii yx=\cou(y)x\quad\forall y\in R.
\end{equation}
Furthermore, the inveribility of $\intgl R$ tells that it is a one dimensional
Yetter-Drinfeld module.

Analogously, the defining equation for left integral elements in $\dad R$
is stated as
$$\lambda\in\intgl{\dad R} \sii \gi\lambda,x\gd 1
	=\gi\lambda,\comm x20\gd\St(\comn x21)\com x1
	\quad\forall x\in R.$$
We note that $\gi\lambda,\comm x20\gd\St(\comn x21)\com x1
=\gi\com{\lambda}0,\com x2\gd\con{\lambda}1\com x1$.
Let now $\lambda\in\intgl{\dad R},\ \lambda\neq 0$. We have an isomorphism
of Yetter--Drinfeld modules $R\simeq\dad R,\ x\mapsto(\lambda\ganid x)$.
Therefore we have the following nondegenerate bilinear form on $R$:
\begin{equation}\label{fornodeg}
(x,y)=(\lambda\ganid x)(y)=\gi\lambda,(\con x1y)\Ss(\com x0)\gd.
\end{equation}

\bigskip

Takeuchi proves also that for a finite braided Hopf algebra the antipode is
an isomorphism as we now sketch. Let $\I=\intgl{\dad R}$. The isomorphism
$\dad R\simeq\I\otimes R$ is given by
$$\alpha=\left(\I\otimes R\to\dad R\otimes R @>\ganid>>\dad R\right).$$
Note that, because of the definition of $\ganid$, $\alpha$ can be
factorized as
$$\I\otimes R \flllad{\id\otimes\Ss} \I\otimes R \flad{\beta} \dad R,$$
which implies that $\I\otimes R \flllad{\id\otimes S} \I\otimes R$ has a left
inverse. Tensoring it with $R$ on the left and composing with the isomorphism
$$R\otimes\I\otimes R \flad{c} \I\otimes R\otimes R\to\dad R\otimes R,$$
we get that $\dad R\otimes R @>\id\otimes\Ss>>\dad R\otimes R$ has a
left inverse. Since $\uni\flad{\dad\cou}\dad R\flad{\dad u}\uni$ is the
identity morphism, we can compose $(\id\otimes\Ss)$ and its left inverse
convenientely with $\dad\cou\otimes\id_R$ and $\dad u\otimes\id_R$, and we
get that $\Ss:R\to R$ has a left inverse.
Since the same argument proves that $\dad{\Ss}:\dad R\to \dad R$ has a
left inverse, this implies that $\Ss$ has a right inverse also.

\section{Braided Hopf algebras of type one} 
\subsection{Semisimplicity of Yetter--Drinfeld categories over
	group algebras}\ \\
\label{loseme}

Let $\Gamma$ be a finite group. Let $H$ be the group
algebra of $\Gamma$ over $\k$, where $\k$ is an algebraically closed
field whose characteristic does not divide the order of $\Gamma$.
We prove that $\yd$ is a semisimple category, and give a complete description
of the simple objects in terms of irreducible representations of some
subgroups of $\Gamma$. This seems to be folklore; it can be found e.g. in
\cite[Prop 3.3]{cibils-rosso:97} in the language of Hopf bimodules
(see also \cite{cibils:97}).
Thanks to \ref{doble}, in order to give all the simple objects of $\yd$
it is enough to give a collection of mutually non isomorphic simple
objects for which the sum of the squares of their dimensions is the
dimension of $\dob H$. It is known, in fact, that the double
of a semisimple and cosemisimple Hopf algebra is semisimple,
see \cite[Cor 10.3.13]{montgomery:93}, but the argument there is not constructive,
in the sense that it refers to Maschke's theorem.

\medskip
We consider the conjugacy classes of $\Gamma$, and choose an element in
each class, which gives a subset $\conor$ of $\Gamma$.
For any $g\in\Gamma$ we denote by $\orb_g=\{xgx^{-1}|x\in\Gamma\}$ the
conjugacy class of $g$, and by $\Gamma_g=\{x\in\Gamma|xg=gx\}$ the isotropy
subgroup of $g$.

\begin{defn}
Let $\rho:\Gamma_g\to\End(V)$ be
an irreducible representation of $\Gamma_g$, and let
$$M(g,\rho):=\Ind_{\Gamma_g}^{\Gamma}V=\k\Gamma\otimes_{\k\Gamma_g}V.$$
For $v\in V,\ x\in\Gamma$, we denote by $\alai vx$ the element
$x\otimes v\in M(g,\rho)$, and by $\alai gx$ the conjugate $xgx^{-1}$.
We take for $M(g,\rho)$ the structure given by
\begin{align*}
& h\ganda \alai vx=\alai v{hx}\quad\mbox{(the induced structure)}, \\
& \delta(\alai vx)=\alai gx\otimes\alai vx,
\end{align*}
which makes $M(g,\rho)$ into an object of $\yd$.
Observe that $\dim M(g,\rho)=[\Gamma:\Gamma_g]\times\dim(\rho)$.
\end{defn}

Given a group $G$ we denote as usually by $\rid G$ the set of
isomorphism classes of irreducible representations of $G$. We often
denote a class in $\rid G$ by a representative element.

\begin{prop}
The objects $M(g,\rho)$ are simple, and any simple object of $\yd$ is
isomorphic to $M(g,\rho)$ for a unique $g\in\conor$ and a unique
$\rho\in\rid{\Gamma_g}$.
\end{prop}

\begin{pf}
Let $g\in\conor,\ \rho\in\rid{\Gamma_g}$. Let $0\neq W\subseteq M(g,\rho)$
be a Yetter--Drinfeld submodule. We have to prove that $W=M(g,\rho)$.
Let $E_g$ be a set of representatives of left coclasses of $\Gamma$ modulo
$\Gamma_g$, i.e. $\Gamma=\bigcup_{x\in E_g}x\Gamma_g$. Observe that
$M(g,\rho)=\bigoplus_{x\in E_g}\k x\otimes V$ as vector spaces, where $V$
is the space affording $\rho$, i.e. $\rho:\Gamma_g\to\Aut(V)$.
Let $0\neq v\in W,\ v=\sum_{x\in E_g}x\otimes v_x=\sum_{x\in E_g}\alai{(v_x)}x$.
Let $p_x=\delta_{\alai gx}\in H^*$ be defined by $p_x(t)=\delta_{\alai gx,t}$.
We have
$$\delta(v)=\sum_{x\in E_g}\alai gx\otimes\alai{(v_x)}x\in H\otimes W
\so\alai{(v_x)}x=(p_x\otimes\id)(\delta(v))\in W\quad\forall x\in E_g.$$
Now, as $v\neq 0$, we have $v_y\neq 0$ for some $y\in E_g$.
Then $v_y=1\otimes v_y=y^{-1}\ganda(\alai{(v_y)}y)\in W$,
but $\k\Gamma_g\ganda v_y=\k 1\otimes V$ because of the ireducibility of
$\rho$, and then $\k 1\otimes V\subseteq W$. Thus
$$\forall x\in E_g,\ \forall v\in V,\quad
\alai vx=(x\ganda v)\in (\Gamma\ganda W)\subseteq W,$$
whence $W=M(g,\rho)$.

\medskip

Let now $h\in\conor$ and $\tau:\Gamma_h\to\End(V')$ be an irreducible
representation of $\Gamma_h$. Define $M(h,\tau)$ as before. If $g\neq h$, it
is immediate that $M(g,\rho)\not\simeq M(h,\tau)$ because $M(g,\rho)$ has
elements of degree $g$ and $M(h,\tau)$ does not. If $g=h$ and
$\rho\not\simeq\tau$ then $M(g,\rho)\not\simeq M(h,\tau)$ because any
isomorphism, being a morphism of comodules, restricts to an isomorphism
between $V$ and $V'$.
Therefore, we have a collection of mutually non isomorphic objects
$M(g,\rho)$ taking $g\in\conor$ and $\rho\in\rid{\Gamma_g}$.
These are all the irreducible objects in $\yd$, since
\begin{align*}
\sum_{g\in\conor}&\sum_{\rho\in\rid{\Gamma_g}}(\dim M(g,\rho))^2
=\sum_{g\in\conor}\sum_{\rho\in\rid{\Gamma_g}}([\Gamma:\Gamma_g]\dim \rho))^2 \\
&=\sum_{g\in\conor}([\Gamma:\Gamma_g]^2\sum_{\rho\in\rid{\Gamma_g}}(\dim \rho)^2)
=\sum_{g\in\conor}([\Gamma:\Gamma_g])^2(\#\Gamma_g) \\
&=\sum_{g\in\conor}(\#\orb_g)^2(\#\Gamma_g)
=\sum_{g\in\conor}(\#\orb_g)(\#\Gamma)
=(\#\Gamma)\sum_{g\in\conor}(\#\orb_g)=(\#\Gamma)(\#\Gamma)=
\dim(\dob H)
\end{align*}
\end{pf}

\begin{cor} \label{abefac}
If $\Gamma$ is abelian, every Yetter--Drinfeld module over $\k\Gamma$
can be decomposed as a direct sum of YD modules of dimension $1$.
\end{cor}
\begin{pf}
It is immediate, since $[\Gamma:\Gamma_g]=\dim\rho=1$ for every $g\in\Gamma$
and $\rho\in\rid{\Gamma_g}$.
\end{pf}

From now on $R$ shall be a Hopf algebra in $\yd$. We shall denote by
$$\prim R=\{x\in R\ |\ \Delta(x)=1\otimes x+x\otimes 1\}$$
the space of primitive elements of $R$.

\begin{lem}
$\prim R$ is a Yetter--Drinfeld submodule of $R$.
\end{lem}

\begin{pf}
Consider the morphism from $R$ to $R\otimes R$ given by
$\id_R\otimes u_R+u_R\otimes\id_R$. It is a morphism in $\yd$, as well as
$\Delta_R$. Then $\prim R$ is the equalizer of both morphisms, which
is a submodule and a subcomodule of $R$.
\end{pf}

We recall from \cite{sweedler:69} that the {\em coradical} of $R$ is
the sum of all its simple subcoalgebras. Since $R$ is in particular a
(classic) coalgebra, we can apply to $R$ the machinery of the coradical
filtration. We shall denote by $R_0$ the coradical of $R$.
We are interested in braided Hopf
algebras $R$ such that $R_0=\k 1$ and $\prim R$ generates $R$ as
an algebra. We give now some first results about such algebras, mainly
for the case in which $\prim R$ is an irreducible object, and postpone to
the next section a more formal and general treatment.

\begin{defn}
Let $R$ be a braided Hopf algebra in $\yd$. We say that $R$ is an
{\em ET-algebra} if $R_0=\k 1$ and $\prim R$ generates $R$.
\end{defn}

\begin{prop} \label{ftabcase}
Let $R$ be an ET-algebra such that $\prim R=M(g,\rho)$ for some
$g\in\Gamma$, $\rho\in\rid{\Gamma_g}$. Then the bosonization $R\# H$
is an extension of the bosonization $R\#\k G$ by the group algebra
$\k(\Gamma/G)$, where $G$ is the subgroup of $\Gamma$ generated by $\orb_g$.
\end{prop}

\begin{pf}
Let $V$ be the space affording $\rho$, i.e. $\rho:\Gamma_g\to\Aut(V)$.
Observe first that $G$ is normal, because if
$h\in\Gamma$ and $g_1,\ldots,g_n\in\orb_g$ then
$\alai {(g_1\cdots g_n)}h=\alai {g_1}h\cdots \alai {g_n}h\in G$.
Observe now that $\delta(R)\subseteq \k G\otimes R$ because $\prim R$
generates $R$, and then $R$ can be considered as a $\k G$-module and a
$\k G$-comodule. Furthermore, it is immediate that $R$ is a Hopf algebra
in $\yds{\k G}$.
Consequently, there exists an inclusion
$$A=R\#\k G\hookrightarrow R\# H=B.$$
Moreover, this inclusion is normal: let $h\in G,\ x\in \Gamma,\ r\in R,\ 
s=\alai vt\in M(g,\rho)\subset R$, where $t\in E_g$ ($=$ left coclasses
$\Gamma/\Gamma_g$). Then
\begin{align*}
\Delta_B(s) &=\Delta_B(s\# 1)=(\id\otimes c\otimes\id)((1\otimes s+s\otimes 1)
	\otimes (1\otimes 1))=\alai gt\otimes s+s\otimes 1, \\
\Ss(s) &=\Ss(s\# 1)=-(\alai gt)^{-1}s, \\
\Delta_B(x) &=\Delta_B(1\# x)=(1\# x)\otimes(1\# x)=
	x\otimes x,
\end{align*}
and thus
\begin{align*}
\Ad_s(h) &=\com s1h\Ss(\com s2)=-{\alai gt}h(\alai gt)^{-1}s+sh=
	-\alai h{(\alai gt)}s+sh\in R\#\k G, \\
\Ad_s(r) &=\com s1r\Ss(\com s2)=-{\alai gt}r(\alai gt)^{-1}s+sr=
	-({\alai gt}\ganda r)s+sr\in R\#\k G, \\
\Ad_x(h) &=\com x1h\Ss(\com x2)=xhx^{-1}=\alai hx\in R\#\k G, \\
\Ad_x(r) &=\com x1r\Ss(\com x2)=xrx^{-1}=x\ganda r\in R\#\k G.
\end{align*}
The condition that $\prim R$ generates $R$ guarantees that
$\forall s\in R,\ \Ad_s(R\# \k G)\subseteq R\#\k G$.
We have therefore a sequence of Hopf algebras
$$\begin{CD}\k @>>> R\#\k G @>>> R\# H @>>> \k(\Gamma/G) @>>> \k.\end{CD}$$
It is straightforward to see that this sequence fulfills the conditions
of \cite[1.2.3]{andrus-devoto:95}, and then the sequence is exact.
\end{pf}

\bigskip

\begin{rem}
The space $\prim R$ may be a simple object in $\yd$, but may be
decomposable when considered as an object in $\yds{\k G}$. For instance,
when $G$ is abelian we know from corollary \ref{abefac} that
$\prim R$ decomposes as a sum of objects of dimension 1.
\end{rem}

\begin{defn}
We say that $R$ {\em can be obtained from the abelian case} if its
bosonization is an extension
$$\begin{CD}\k @>>> R\#\k\Gamma_1 @>>> R\# H @>>> \k\Gamma_2 @>>> \k\end{CD}$$
where $\Gamma_1$ is abelian.
\end{defn}

\begin{lem}
Let $R$ be an ET-algebra such that $\prim R=M(g,\rho)$ for some $g$,
and $\rho\in\rid{\Gamma_g}$. Let $G$ be the subgroup generated by $\orb_g$.
If $\Gamma_g\triangleleft\Gamma$ then $G$ is abelian, and thus $R$ can be
obtained from the abelian case.
\end{lem}

\begin{pf}
Let $t\in\orb_g$, $t=\alai gx$. Then we have $\Gamma_t=\alai {\Gamma_g}x
=x\Gamma_gx^{-1}=\Gamma_g$, which implies that $\Gamma_{t_1}=\Gamma_{t_2}$
for any $t_1,t_2\in\orb_g$. Thus any two elements in $\orb_g$ commute,
and hence the group generated by $\orb_g$ is abelian.
\end{pf}

\begin{exmp}
The preceding lemma has the following application: if all the
isotropy subgroups of $\Gamma$ are normal, then any ET-algebra
with an irreducible space of primitive elements can be obtained
from the abelian case. This happens, for instance, for $\DD_4$.
Other examples are the groups such that every subgroup is normal.
It is known that these groups are abelian, or
of the form $H\times A$, where $H$ is the quaternion group, i.e.
$$H=\{1,-1,i,-i,j,-j,k,-k|\ i^2=j^2=k^2=-1,\ ij=k,\ jk=i,\ ki=j\},$$
and $A$ is an abelian group without elements of order $4$
(see \cite{carmichael:56}).
\end{exmp}

\begin{prop}
Let $R$ be an ET-algebra such that $\dim\prim R=2$. Then $R$ can be obtained
from the abelian case.
\end{prop}

\begin{pf}
Let $M=\prim R$. We have $\dim M=2$ and then there are three possibilities:
\begin{enumerate}
\item $M$ is decomposable as $M=M(g_1,\chi_1)\oplus M(g_2,\chi_2)$ with $g_i$
	central in $\Gamma$ and $\chi_i$ characters.
\item $M$ is simple, and then $M=M(g,\rho)$ with $\Gamma_g=\Gamma$ and
	$\dim\rho=2$, or
\item $M=M(g,\chi)$ with $[\Gamma:\Gamma_g]=2$ and $\chi$ a character.
\end{enumerate}

In the first case, let $G$ be the group generated by $g_1$ and $g_2$. It is
abelian because the $g_i$ are central. The construction of proposition
\ref{ftabcase} can be made with this $G$.

In cases 2 and 3 we have $\Gamma_g=\Gamma$ or $[\Gamma:\Gamma_g]=2$,
and the result follows from the lemma above.
\end{pf}

\subsection{Bialgebras of type one}\ \\
As we said before, we are interested in certain classes of braided
Hopf algebras, which we define in this section. Most of the following
can be made in any braided category, as it is done by Schauenburg in
\cite{schauenburg:96}. Given a braided category $\cate$, he is obliged to work
in an $\NN$-graded category $\cate^\NN$. To avoid these technicalities
we shall work in $\yd$ and $\ydi$.

\begin{defn}\label{defgradbi}
A graded Hopf algebra in $\yd$ or $\ydi$ is simply a Hopf algebra $R$
in any of these categories such that $R=\bigoplus_n\agg Rn$ and
$$\agg Ri\agg Rj\subseteq\agg R{i+j},\quad
\Delta(\agg Rk)\subseteq\bigoplus_{i+j=k}\agg Ri\otimes\agg Rj.$$
\end{defn}

\bigskip

An important application of the existence of an integral is the following
\begin{prop}[Nichols] \label{dualpoin}
Let $R=\oplus_{i=0}^N\agg Ri$ be a graded Hopf algebra in $\yd$ (it is
in particular finite dimensional), and suppose that $\agg RN\neq 0$. Then
$\dim\agg Ri=\dim\agg R{N-i}\ \forall i$.
\end{prop}
\begin{pf}
Since $R$ is graded, it is clear that $\dad R$ is also a graded Hopf algebra
in $\yd$. Let $\lambda\in\dad R$ be a non zero left integral. We have then
$\lambda=\sum_i\lambda_i$, where $\lambda_i\in\agg{\dad R}i$ is the component 
of degree $i$. It is immediate, looking at \eqref{rinte}, that each
$\lambda_i$ is a left integral in $\dad R$. Then, by the one dimensionality
of the space of integrals, we have $\lambda=\lambda_j$ for some $j$.
We recall now (see \eqref{fornodeg}) that $\lambda$ defines a non degenerate
bilinear form
$$(x,y)=\gi\lambda,(\con x1y)\Ss(\com x0)\gd.$$
Since $\agg Ri$ and $\agg Rk$ are orthogonal if $i+k=j$, this map
induces a non degenerate bilinear form between $\agg Ri$ and
$\agg R{j-i}$ for each $i$. Hence in particular we have that
$\agg Ri=0\ \forall i>j$, and then $j=N$, whence the thesis.
\end{pf}

\begin{defn}\label{deftoba}
A {\em braided Hopf algebra of type one}, or briefly {\em TOBA},
\footnote{The Tobas are an aboriginal ethnic group living in the north
of Argentina.}
in $\yd$ or in $\ydi$ is a graded Hopf algebra in any of
these categories such that
\begin{align}
\label{toba1} \agg R0 &\simeq\k, \\
\label{toba2} \left(\oplus_{i\ge 1}\agg Ri\right)^2 &= 
			\oplus_{i\ge 2}\agg Ri, \\
\label{toba3} \prim R &= \agg R1.
\end{align}
\end{defn}

\begin{rem} \label{grbiesho}
A graded bialgebra in $\yd$ or $\ydi$ which satisfies these conditions
is automatically a Hopf algebra, thanks to \cite[Lemma 5.2.10]{montgomery:93}.
\end{rem}

It is easy to see that if $R$ is a TOBA then the unit and counit are
respectively the canonical inclusion and canonical projection
$$u:\k=\agg R0\hookrightarrow R,\quad\cou:R\sobre\agg R0=\k.$$
It is easy to see that in presence of \eqref{toba1} the
condition \eqref{toba3} is equivalent to the condition
$$R_1=R_0\wedge R_0=\agg R0\oplus\agg R1,$$
where $\wedge$ stands for the wedge product and $R_0\subset R_1\subset\dots$
stands for the coradical filtration of $R$ (see \cite[Ch. 9]{sweedler:69}).
Moreover, it is easy to see by induction that the condition \eqref{toba2}
is equivalent to
$$(\agg R1)^n=\agg Rn\quad\forall n\ge 1,$$
which in presence of \eqref{toba3} can be stated by saying that
$\prim R$ generates $R$.

\begin{exmp}
Let $R=\k[x]/x^{(p^2)}$, where $\car\k=p$. The comultiplication
is determined by imposing $x$ to be a primitive element.
This is a (usual) graded Hopf algebra which verifies \eqref{toba1} and
\eqref{toba2}, but not \eqref{toba3}. Its dual is a graded Hopf algebra
which verifies \eqref{toba1} and \eqref{toba3} but not \eqref{toba2}.
Another example is the tensor algebra $T(V)$ of a vector space $V$ of
dimension greater than $1$, with comultiplication determined by
$V\subseteq \prim{T(V)}$. This Hopf algebra satisfies \eqref{toba1}
and \eqref{toba2} but not \eqref{toba3}. Indeed, $\prim{T(V)}$ is
the free Lie algebra generated by $V$.

It is not known whether a finite dimensional graded (braided) Hopf
algebra satisfying \eqref{toba1} and \eqref{toba3} should satisfy
\eqref{toba2} provided that $\car\k=0$.
This was proved in \cite{andrus-schneider:97} in the case $\dim\prim R=1$.
\end{exmp}

\bigskip

We give three ways to construct a TOBA. The second one is due
to Nichols (see \cite{nichols:78}), from where we borrow the name.
The first one can be seen as a rewriting of that of Nichols in the language
of braided categories, and is due to Schauenburg (see \cite{schauenburg:96},
see also \cite{rosso:95,rosso:92}, \cite{bespa-draba:96}).
The last one is inspired in the work of
Lusztig \cite{lusztig:93} and is stated for the category $\ydi$, where
$H=\k\Gamma$ (Lusztig's algebras $\f$ and ${}'\f$ are braided Hopf
algebras in a category of comodules).
The approach of \cite{schauenburg:96} is in fact motivated by this work.
We prefer to give the way of \cite{schauenburg:96} first because it seems
more useful to us, and then give that of \cite{nichols:78} in the terms
of \cite{schauenburg:96}.
It is important to note that we work in $\yd$ and $\ydi$,
rather than in an $\NN$-graded category.

\bigskip

Let $n\in\NN,\ n\ge 2$. Let $\SS_n$ and $\BB_n$ be the symmetric and
braid groups defined in \ref{defbn}.
There is a projection $\BB_n\sobre\SS_n,\ \sigma_i\mapsto\tau_i$.

Let $x\in\SS_n$. We denote by $\ell(x)$ the length of a minimal word
in the alphabet $\{\tau_i|\ 1\le i< n\}$ which represents $x$.
For $y\in\BB_n$ we denote also by $\ell(y)$ the length of a minimal
word in the alphabet $\{\sigma_i,\ \sigma_i^{-1}|\ 1\le i<n\}$ which
represents $y$. There is a unique section $s:\SS_n\to\BB_n$
to the projection $\BB_n\to\SS_n$ such that
$s(\tau_i)=\sigma_i$ and $s(ww')=s(w)s(w')$ whenever
$\ell(w\cdot w')=\ell(w)+\ell(w')$. It is given by
\begin{equation} \label{defsec}
(w=\tau_{i_1}\cdots\tau_{i_j})\mapsto(\sigma_{i_1}\cdots\sigma_{i_j})
	\quad\mbox{if }\ell(w)=j\qquad\mbox{(thus $\ell s=\ell$)}.
\end{equation}
It is clear that it is unique; it is proved in \cite[64.20]{curtis-reiner:94}
that it is well defined. Using this section, we define the $\SSS$ morphisms:
let $V$ be an object in $\yd$. As in remark \ref{bnactua}, $\BB_n$ acts
on $V^{\otimes n}$. For $w\in\BB_n$ we denote also by $w$ the
corresponding morphism given by this action.
If $X\subseteq\SS_n$, we define the morphism
$$\SSS_X:V^{\otimes n}\to V^{\otimes n},\quad \SSS_X=\sum_{x\in X}s(x).$$
Let $k_1,\ldots,k_j\in\NN$ such that $k_1+\dots+k_j=n$.
We denote by $X_{k_1,\ldots,k_j}\subseteq\SS_n$ the
$(k_1,\ldots,k_j)$-shuffle and $Y_{k_1,\ldots,k_j}\subseteq\SS_n$
the inverse of $X_{k_1,\ldots,k_j}$, i.e.
\begin{align*}
X_{k_1,\ldots,k_j}&=\{x\in\SS_n|\,x^{-1}(k_1+\dots+k_i+1)<\dots
	<x^{-1}(k_1+\dots+k_{i+1})\ \forall i=0,\ldots,j-1\} \\
Y_{k_1,\ldots,k_j}&=X_{k_1,\ldots,k_j}^{-1}=
\{x^{-1}|\,x\in X_{k_1,\ldots,k_j}\}
\end{align*}
We then define $\SSS_{k_1,\ldots,k_j}=\SSS_{X_{k_1,\ldots,k_j}}$,
$\SSS^n=\SSS_{1,1,\ldots,1}=\SSS_{\SS_n}$, and 
$\TTT_{k_1,\ldots,k_j}=\SSS_{Y_{k_1,\ldots,k_j}}$.

Then, for instance, $\SSS^2=\id+c$, and
$$\SSS_{2,1}=\id_{V^{\otimes 3}}+\id_V\otimes c_{V,V}+(\id_V\otimes c_{V,V})
	(c_{V,V}\otimes\id_V):V^{\otimes 3}\to V^{\otimes 2}\otimes V.$$
We observe that for $i+j=n$, $\SSS^n=(\SSS^i\otimes\SSS^j)\SSS_{i,j}$.

\begin{defn}
Let $V$ be an object of $\yd$. We denote by $T^n(V)=V^{\otimes n}$ and
by $T(V)$ the tensor object
$$T(V)=\k\oplus V\oplus V^{\otimes 2}\oplus\cdots\oplus
	V^{\otimes n}\oplus\cdots$$
$T(V)$ is not an object of $\yd$, but an object of $\ydi$.
We consider on $T(V)$ two different bialgebra structures, which we
denote by $A(V)$ and $C(V)$. Both are graded bialgebras in the sense
of \ref{defgradbi}, and we denote
$$m=\mathop{\oplus}_{i,j}m_{i,j},\quad
	\Delta=\mathop{\oplus}_{i,j}\Delta_{i,j},$$
where $m_{i,j}:A^i(V)\otimes A^j(V)\to A^{i+j}(V)$ and
$\Delta_{i,j}:A^{i+j}(V)\to A^i(V)\otimes A^j(V)$ and the same for $C(V)$.
We take for both $A(V)$ and $C(V)$ the unit and counit given by
inclusion $\k\to T(V)$ and projection $T(V)\to\k$.
We take on $A(V)$ the multiplication given by
$$m_{i,j}=\id:A^i(V)\otimes A^j(V)\to A^{i+j}(V).$$

\medskip
\noi There exists only one comultiplication making $A(V)$ into a bialgebra
in $\ydi$ with
$$\Delta_{1,0}=\id:V\to V\otimes\k,\quad
\Delta_{0,1}=\id:V\to\k\otimes V,$$
which is given by
$$\Delta_{i,j}=\SSS_{i,j}:A^{i+j}(V)\to A^i(V)\otimes A^j(V).$$
Dually, we take on $C(V)$ the comultiplication given by
$$\Delta_{i,j}=\id:C^{i+j}(V)\to C^i(V)\otimes C^j(V).$$

\medskip
\noi There exists only one multiplication making $C(V)$ into a bialgebra
in $\ydi$ with
$$m_{0,1}=\id:\k\otimes V\to V,\ m_{1,0}=\id:V\otimes\k\to V,$$
which is given by
$$m_{i,j}=\TTT_{i,j}:C^i(V)\otimes C^j(V)\to C^{i+j}(V).$$
\end{defn}

\medskip
\noi There exists only one morphism of bialgebras $\SSS:A(V)\to C(V)$
such that $\SSS|_{A^1}=\id:V\to V$. This is the graded morphism given by
$$\SSS=\oplus_n\SSS^n:A^n(V)\to C^n(V).$$
Let $B(V)=\bigoplus_nB^n(V)$ be the image $\SSS(A(V))\subset C(V)$.
This is a bialgebra in $\ydi$ and is isomorphic to
$$A(V)/\ker(\SSS)=\bigoplus_n(A^n(V)/\ker(\SSS^n)).$$

We have then a graded bialgebra $B(V)$ in $\ydi$ which, by construction,
verifies \eqref{toba1} and \eqref{toba2}. It also verifies \eqref{toba3}
since, as a subbialgebra of $C(V)$, the comultiplication components
$\Delta_{i,j}:B^{i+j}\to B^i\otimes B^j$ are injective for all $i,j\in\NN$.
Hence we have

\medskip
\begin{defn}
Let $V\in\yd$, $\toba V:=\bigoplus_n(A^n(V)/\ker\SSS^n)\subseteq C(V)$.
It is a TOBA with $\prim{\toba V}\simeq V$. As we noted in \ref{grbiesho}
it has an antipode. It is given by $\Ss(x)=-x\quad\forall x\in\agg R1$
and it is extended by \ref{sdavuel}. The following proposition proves
that a TOBA is fully determined by its space of primitive elements,
and thus $\toba V$ can be defined alternatively by conditions
\eqref{toba1}--\eqref{toba3} plus $\agg{\toba V}1=V$.
\end{defn}

\bigskip
\begin{prop} \label{prfully}
Let $R$ be a TOBA. Let $V=\agg R1$, $p:A(V)\to R$ be the algebra
surjection induced by the inclusion $V\hookrightarrow R$, and $I$ be
its kernel. Then $I=\ker(\SSS)$.
\end{prop}

\begin{pf}
We prove first that $I\supseteq\ker(\SSS)$.
Since both $I$ and $\ker(\SSS)$ are homogeneous, we have to prove that
$I_n\supseteq\ker(\SSS^n)$, where $I_n$ is the homogeneous component of
$I$ of degree $n$. We proceed by induction.

For $n=1$ there is nothing to prove, since $\SSS^1=\id$.
Let $p:A(V)\to R$ be the projection, and suppose that the inclusion
is true for $m<n$. Let $x\in\ker(\SSS^n)$. We have that
$$\Delta(x)=\sum_{k+l=n}\SSS_{k,l}x\in\sum_{k+l=n}A^k\otimes A^l,$$
but $\SSS^n=(\SSS^k\otimes\SSS^l)\SSS_{k,l}$, and hence
$$\SSS_{k,l}(x)\in\ker(\SSS^k\otimes\SSS^l)=
\ker\SSS^k\otimes V^{\otimes l}+V^{\otimes k}\otimes\ker\SSS^l,$$
whence $\SSS_{k,l}(x)\in I\otimes A+A\otimes I$ if $k,l<n$.
Then $(p\otimes p)(\SSS_{k,l}(x))=0$ if $k,l<n$, which implies that
$$\Delta(p(x))=\sum_{k+l=n}(p\otimes p)\SSS_{k,l}(x)=
(p\otimes p)(\SSS_{n,0}(x)+\SSS_{0,n}(x))=p(x)\otimes 1+1\otimes p(x).$$
Thus $p(x)\in\prim R$, but $p(x)\in\agg Rn$ and $n>1$, whence
$p(x)=0$ and then $x\in I$. This proves the first inclusion.

We have now the quotient morphism of coalgebras (of braided Hopf algebras,
in fact)
$$A(V)/\ker\SSS \to R,$$
which is injective on $(A(V)/\ker(\SSS))_1$ (the second term of the
coradical filtration), since
$$(A(V)/\ker(\SSS))_1=\k\oplus V=\agg R0\oplus \agg R1.$$
By \cite[5.3.1]{montgomery:93}, the quotient morphism is injective, which
says that $\ker(\SSS)=I$.
\end{pf}

\begin{rem}\label{Rdepdec}
We note that $R=\toba V$ depends as a braided Hopf algebra only on the
braiding $c_{V,V}\in\End(V\otimes V)$. This allows to consider $R$ in
different categories, as long as $c_{V,V}$ remains unchanged.
\end{rem}

We give now a second construction of a TOBA. See \cite{nichols:78} for details.
\begin{defn}[Nichols]
A graded bialgebra $B=\bigoplus_{i\ge 0}\agg Bi$ is called {\em a bialgebra
of type one} if it verifies the following conditions:
\begin{align}
\label{toa2} \left(\oplus_{i\ge 1}\agg Bi\right)^2 &= \oplus_{i\ge 2}\agg Bi, \\
\label{toa3} \agg B0\wedge \agg B0 &=\agg B0\oplus \agg B1.
\end{align}
We define similarly the notion of {\em Hopf algebra of type one}.
\end{defn}

\begin{rem}
Let $B=\bigoplus_{i\ge 0}\agg Bi$ be a bialgebra. As in the braided
case, it follows from \cite[Lemma 5.2.10]{montgomery:93} that $B$
has an antipode if and only if $\agg B0$ does.
\end{rem}

Nichols constructs bialgebras of type one starting out from Hopf
bimodules. We relate his construction to that of
Schauenburg. In order to do this we need the following

\begin{lem}
Let $H$ be a Hopf algebra, $R=\bigoplus_{n\ge 0}\agg Rn$
be a graded Hopf algebra in $\yd$, and $A=R\# H$.
This is a graded Hopf algebra with respect to the grading
$$A=\bigoplus_{n\ge 0}\agg An,\ \agg An=(\agg Rn\otimes H)\subseteq R\# H.$$
If $R$ is a TOBA then $A$ is a bialgebra of type one such that
$\agg A0\simeq H$.
Conversely, let $B=\bigoplus_{n\ge 0}\agg Bn$ be a graded Hopf algebra.
We have the canonical morphisms of Hopf algebras $\agg B0\hookrightarrow B$
and $B\sobre\agg B0$. Let $R=\coind B{\pi}$; it is a graded Hopf algebra
in $\yds {\agg B0}$. Hence, if $B$ is a Hopf algebra of type one
then $R$ is a TOBA.
\end{lem}

\begin{pf}
Condition \eqref{toba1} is easily seen to be equivalent to the condition
$\agg A0=H$. Then the equivalence between \eqref{toba2} and \eqref{toa2}
is a consequence of the following: let $M$ and $N$ be subspaces of $R$.
We claim that if $N$ is an $H$-submodule then $(MN)\# H=(M\# H)(N\# H)$.
For this, let $m\in M,\ n\in N,\ h\in H$.
Then $(mn\# h)=(m\#1)(n\# h)\in (M\# H)(N\# H)$, which implies one inclusion.
The other is immediate under the hypothesis of $N$ being a submodule.

As we remarked after the definition \ref{deftoba}, it is easy to see that
\eqref{toba3} is equivalent to the condition $\k\wedge\k=\agg R0\oplus\agg R1$.
The equivalence between conditions \eqref{toba3} and \eqref{toa3} is a
consequence of the following: let $M$ and $N$ be subspaces of $R$.
We claim that if $N$ is an $H$-subcomodule then
$(M\# H)\wedge (N\# H)=(M\wedge N)\# H$. To see this, we consider both
subspaces as kernels of certain morphism: let us denote by $\tau$
the usual flip $x\otimes y\mapsto y\otimes x$. If $X$ is a subspace
of $Y$, we denote by $p_X:Y\to Y/X$ the canonical projection. Since
$N$ is a subcomodule, $R/N$ has an $H$-comodule structure, which we
denote also by $\delta$. Thus,
\begin{align*}
(M\wedge &N)\# H = \ker((p_M\oti p_N)\Delta_R) \# H \\
&= \ker((p_M\oti p_N\oti\id)(\Delta_R\oti\id)) \\
&= \ker((p_M\oti\id\oti p_N\oti\id)(\id\oti\tau\oti\id)
	(\Delta_R\oti\Delta_H)) \\
&= \ker((\id\oti m_H\oti\id\oti\id)(\id\oti\tau\oti\id\oti\id)
	(\id\oti\id\oti\delta\oti\id)(p_M\oti\id\oti p_N\oti\id)
	(\id\oti\tau\oti\id)(\Delta_R\oti\Delta_H)) \\
&= \ker((p_M\oti m_H\oti p_N\oti\id)(\id\oti\id\oti\tau\oti\id)
	(\id\oti\delta\oti\id\oti\id)(\Delta_R\oti\Delta_H)) \\
&= \ker((p_M\oti\id\oti p_N\oti\id)(\id\oti c\oti\id)
	(\Delta_R\oti\Delta_H)) \\
&= \ker((p_{(M\# H)}\oti p_{(N\# H)})\Delta_{R\# H}) \\
&= (M\# H)\wedge (N\# H).
\end{align*}
\end{pf}

Let $(P,\delta_P:P\to P\otimes H)$ be a right $H$-comodule, and
$(Q,\delta_Q:Q\to H\otimes Q)$ be a left $H$-comodule. We denote as
usual by $\coti_H$ the cotensor product, i.e.
$$P\coti_HQ=\equ(P\otimes Q 
\mathop{\rightrightarrows}^{\delta_P\otimes\id}_{\id\otimes\delta_Q}
P\otimes H\otimes Q)=\{\sum p_i\otimes q_i | \sum \com{(p_i)}0\otimes
\com{(p_i)}1\otimes q_i=p_i\otimes\con{(q_i)}1\otimes \com{(q_i)}0\}. $$

Let $H$ be a Hopf algebra and $M\in\hobim$ (see \ref{yd-hobim}).
We denote by
\begin{align*}
A_H(M)&=T_H(M)=H\oplus M\oplus (M\otimes_H M)\oplus (M\otimes_HM\otimes_HM)
\oplus\cdots=\bigoplus_{i\ge 0}A_H^i(M), \\
C_H(M)&=T_H(M)=H\oplus M\oplus (M\coti_H M)\oplus (M\coti_HM\coti_HM)
\oplus\cdots=\bigoplus_{i\ge 0}C_H^i(M).
\end{align*}
As before, $A_H(M)$ (resp. $C_H(M)$) has a canonical graded multiplication
(resp. comultiplication) given by projection
$$A_H^i(M)\otimes A_H^j(M)\to A_H^{i+j}(M)=A_H^i(M)\otimes_HA_H^j(M)$$
(resp. inclusion $C_H^{i+j}=C_H^i\coti_H C_H^j\to C_H^i\otimes C_H^j$).
Moreover, $A_H(M)$ can be endowed with a (unique) comultiplication
which makes it into a bialgebra such that in degree $1$ it is given by
$$A_H^1(M)=M @>\delta_l+\delta_r>> (H\otimes M)\oplus (M\otimes H)=
	\agg {[A_H(M)\otimes A_H(M)]}1,$$
and $C_H(M)$ can be endowed with a (unique) multiplication
which makes it into a bialgebra such that in degree $1$ it is given by
$$\agg {[C_H(M)\otimes C_H(M)]}1=(H\otimes M)\oplus (M\otimes H)
	@>m_l+m_r>> M=C_H^1(M).$$

There exists a unique bialgebra map $A_H(M)\to C_H(M)$ which is the
identity on degrees $0$ and $1$. We denote by $B_H(M)$ its image.
This is a bialgebra of type one. Moreover, since $H$ is a Hopf algebra,
$B_H(M)$ is a Hopf algebra.
This construction is related to that of Schauenburg by the following diagram.
For $\cate$ a braided category, we denote by $\hopf{\cate}$ the subcategory
of the Hopf algebras in $\cate$, by $\#$ the bosonization functor and by $S$
the functor $\yd\to\hobim$ of proposition \ref{yd-hobim}. In this diagram we
denote also by $B$ (instead of ${\frak t}$) the functor giving the
TOBA in $\yd$. Then we have a commutative diagram
$$\begin{CD}
\yd @>B>> \hopf{\yd} \\
@VSVV @V\#VV \\
\hobim @>B>> \hopf{\modui \k}.
\end{CD}$$
The proof that this diagram commutes is straightforward but tedious.
One can verify that the diagram commutes replacing $B$ with $A$ and
$C$, the tensor and cotensor bialgebras, and then note that
the following diagram commutes $\forall V\in\yd$:
$$\begin{CD}
(AV)\# H @>>> A(SV) \\
@VVV @VVV \\
(CV)\# H  @>>> C(SV),
\end{CD}$$
where the left and right sides are the (universal) morphisms $A\to C$,
and the top and bottom sides are the natural equivalences given by the
commutativity of the first diagram with $B$ replaced by $A$ and $C$
respectively.

\bigskip
\begin{rem}
Let $V$ be a $\k$-vector space and $c\in\Aut(V\otimes V)$,
satisfying the braid equation, namely
$$(c \otimes \id)(\id \otimes c)(c \otimes \id)
	=(\id \otimes c)(c \otimes \id)(\id \otimes c).$$
We remark that we can define $\toba V=AV/\ker(\SSS)$ in the same vein
as before, where $\SSS(v)=\sum_{x\in\SS_n}s(x)(v)$ for $v\in A^n(V)$,
and the group $\BB_n$ acting on $V^{\otimes n}$ via $c$.
\end{rem}

\medskip
The last way we give to construct a TOBA is by means of a bilinear form on
$A(M)$. The idea is the same Lusztig uses to construct the algebra $\f$,
and is in fact the motivation for Schauenburg to construct the morphism
$\SSS$ (see \cite{lusztig:93}, \cite{schauenburg:96}). M\"uller uses this
presentation to prove that the nilpotent part $\n_+$ of the
Frobenius-Lusztig kernel $\u$ is a TOBA over $\u_0$ (see \cite{muller:98}).
In our context the form is not a pairing between $A(V)$ and itself, but
between $A(V)$ and $A(W)$, $W$ being another (possibly the same) vector
space. We begin with a useful result.

\medskip
\renewcommand{\theenumi}{(\alph{enumi})}

\begin{lem}\label{SnSn}
Let $U,Z$ be $\k$-vector spaces with an action of $\BB_n$ (in the usual
cases $U=V^{\otimes n}$, $Z=W^{\otimes n}$). We denote for $u\in U$
$$\SSS^nu=\sum_{x\in\SS_n}s(x)(u),$$
and analogously for $z\in Z$. Suppose we have a bilinear form
$(|):U\otimes Z\to\k$ such that either
\begin{enumerate}
\item \label{casoa} $(\sigma_i(u)|z)=(u|\sigma_i(z))$ or
\item \label{casob} $(\sigma_i(u)|z)=(u|\sigma_{n-i}(z))$.
\end{enumerate}
Then we have $(\SSS^nu|z)=(u|\SSS^nz)$ for $u\in U$, $z\in Z$.
\end{lem}

\begin{pf}
Let $x\in\SS_n$, $x=\tau_{i_1}\cdots\tau_{i_d}$ with $\ell(x)=d$.
Then $s(x)=\sigma_{i_1}\cdots\sigma_{i_d}$ and
$s(x^{-1})=\sigma_{i_d}\cdots\sigma_{i_1}$, since
$x^{-1}=\tau_{i_d}\cdots\tau_{i_1}$ and $\ell(x^{-1})=\ell(x)=d$.
Furthermore, let $T\in\SS_n$ be defined by \linebreak
$T:\{1,\ldots,n\}\to\{1,\ldots,n\}$, $T(i)=n+1-i$. Let $D$ be the inner
automorphism defined by $T$, that is, $D:\SS_n\to\SS_n,\ D(x)=TxT^{-1}$.
We observe that $D(\tau_i)=\tau_{n-i}$ for $1\le i\le n-1$. Moreover, since
$T^2=\id$, we have $D^2=\id$. Thus, if $x\in\SS_n$,
$x=\tau_{i_1}\cdots\tau_{i_d}$ and $\ell(x)=d$, we have
$D(x^{-1})=D(\tau_{i_d}\cdots\tau_{i_1})=\tau_{n-i_d}\cdots\tau_{n-i_1}$,
whence $\ell(D(x^{-1}))\le\ell(x)$.
Since $D((D(x^{-1}))^{-1})=(D^2(x^{-1}))^{-1}=x$ and the previous inequality
holds true $\forall x\in\SS_n$, we have $\ell(D(x^{-1}))=\ell(x)$.
Thus,
$$s(D(x^{-1}))=s(\tau_{n-i_d}\cdots\tau_{n-i_1})
	=\sigma_{n-i_d}\cdots\sigma_{n-i_1}.$$
We have now in case \ref{casoa} that
$$(s(x)(u)|z)=(\sigma_{i_1}\cdots\sigma_{i_d}(u)|z)
	=(u|\sigma_{i_d}\cdots\sigma_{i_1}(z))=(u|s(x^{-1})(z)),$$
and hence 
$$(\SSS^nu|z)=\sum_{x\in\SS_n}(s(x)u|z)
	=\sum_{x\in\SS_n}(u|s(x^{-1})z)=(u|\SSS^nz).$$
In case \ref{casob}, we have
$$(s(x)(u)|z)=(\sigma_{i_1}\cdots\sigma_{i_d}(u)|z)
	=(u|\sigma_{n-i_d}\cdots\sigma_{n-i_1}(z))=(u|s(D(x^{-1}))(z)),$$
and hence
$$(\SSS^nu|z)=\sum_{x\in\SS_n}(s(x)u|z)
	=\sum_{x\in\SS_n}(u|s(D(x^{-1}))z)=(u|\SSS^nz).$$
\end{pf}

\begin{rem} \label{bfhe}
Bilinear forms as in \ref{SnSn} happen to exist very often.
Suppose for instance that we have a bilinear form $(|):V\otimes W\to\k$
such that one of the following cases arises:
\begin{enumerate}
\item $(c_V(v_1\otimes v_2)|w_1\otimes w_2)
	=(v_1\otimes v_2|c_W(w_1\otimes w_2))$ for the form
	$(v_1\otimes v_2|w_1\otimes w_2)=(v_1|w_1)(v_2|w_2)$.
\item $(c_V(v_1\otimes v_2)|w_1\otimes w_2)
	=(v_1\otimes v_2|c_W(w_1\otimes w_2))$ for the form
	$(v_1\otimes v_2|w_1\otimes w_2)=(v_1|w_2)(v_2|w_1)$.
\end{enumerate}
In case \ref{casoa} we define
$$(\nftp v1n|\nftp w1n)_>=\prod_i(v_i|w_i),$$
and this fits into case \ref{casoa} of \ref{SnSn}.

In case \ref{casob} we define
$$(\nftp v1n|\nftp w1n)_<=\prod_i(v_i|w_{n+1-i}),$$
and this fits into case \ref{casob} of \ref{SnSn}.

It is clear that if $(|)$ is non degenerate, then $(|)_>$
(resp. $(|)_<$) is non degenerate.
These cases are satisfied in the following examples.
\end{rem}

\begin{exmp} \label{ejelus}
Let $(i\cdot j)_{1\le i,j\le n}$ be a Cartan datum (see \cite{lusztig:93}
for the definition), and let $q$ be an indeterminate over $\k$.
We take $V=\k(q)\theta_1\oplus\ldots\oplus\k(q)\theta_n$, and
define $c_V(\theta_i\otimes\theta_j)=q^{i\cdot j}\theta_j\otimes\theta_i$.
Furthermore, we take $(|):V\otimes V\to\k(q)$ given by
$$(\theta_i|\theta_j)=(1-q^{-2i\cdot i})^{-1}\delta_{i,j}.$$
It is easy to see that this is a non degenerate bilinear form such that
$$(c(\theta_{i_1}\otimes \theta_{i_2})|\theta_{j_1}\otimes \theta_{j_2})
=(\theta_{i_1}\otimes \theta_{i_2}|c(\theta_{j_1}\otimes \theta_{j_2})),$$
whence we are in case \ref{casoa} of the above remark.
\end{exmp}

\begin{exmp}
Let $\cate$ be a braided abelian rigid
category which can be embedded in $\cate'$, a braided abelian category
in which countable direct sums exist. This is the case for instance of
$\yd\hookrightarrow\ydi$ for $H$ any Hopf algebra. Let $\dai V$ be the
left dual of $V$ in $\cate$, and $(|):V\otimes\dai V\to\k$ be the
evaluation map. Lemma \ref{cphi} tells that this fits into case \ref{casob}
of the remark.
\end{exmp}

\renewcommand{\theenumi}{\arabic{enumi}}
\medskip
\begin{defn} \label{defcor}
Let $U,Z$ be $\k\BB_n$-modules with a bilinear form $(|):U\otimes Z\to\k$.
We denote
$$[,]:U\otimes Z\to\k,\quad [u,z]=(\SSS^nu|z).$$
According to this, for $V,W$ $\k$-vector spaces with braidings
$c_V,c_W$ and a bilinear form $(|):V\otimes W\to\k$ satisfying
\ref{casoa} of remark \ref{bfhe} (resp. \ref{casob}), we define
$[,]:AV\otimes AW\to\k$ by
\begin{enumerate}
\item $[1,1]=1$.
\item $[u,z]=0$ if $u\in A^iV$, $z\in A^jW$ and $i\neq j$.
\item $\left\{\begin{array}{ll}
	[u,z]=[u,z]_>=(\SSS^nu|z)_> & \mbox{ if } u\in A^nV,\ z\in A^nW
		\mbox{ for the case \ref{casoa}} \\
	{}[v,w]=[v,w]_<=(\SSS^nv|w)_< & \mbox{ if } u\in A^nV,\ z\in A^nW
		\mbox{ for the case \ref{casob}}
	\end{array}\right.$
\end{enumerate}
\end{defn}

\begin{lem} \label{relpc}
Let $V,W$ be as above. Let us suppose that we are in case \ref{casoa}
(resp. \ref{casob}) of remark \ref{bfhe}. Then we have respectively
\begin{align*}
\mbox{\ref{casoa} } &[u,z\cdot z']_>=[\com u1,z]_>[\com u2,z']_>,\quad
	[u\cdot u',z]_>=[u,\com z1]_>[u',\com z2]_>. \\
\mbox{\ref{casob} } &[u,z\cdot z']_<=[\com u1,z']_<[\com u2,z]_<,\quad
	[u\cdot u',z]_<=[u,\com z2]_<[u',\com z1]_<.
\end{align*}
\end{lem}

\begin{pf}
For $u\in A^nV$ and $i+j=n$, we denote
$$(\SSS_{i,j}(u))_i\otimes(\SSS_{i,j}(u))_j=\SSS_{i,j}(u)
	\in A^iV\otimes A^jV.$$
In case \ref{casoa} we have, for $z\in A^iW$, $z'\in A^jW$ and
$u\in A^nV$,
\begin{align*}
[u,z\cdot z'] &= (\SSS^nu|z\cdot z')
	=((\SSS^i\otimes\SSS^j)(\SSS_{i,j}u)|z\cdot z') \\
	&=(\SSS^i(\SSS_{i,j}u)_i|z)(\SSS^j(\SSS_{i,j}u)_j|z').
\end{align*}
From the other hand, we have
\begin{align*}
[\com u1,z][\com u2,z']&=\sum_{k+l=n}[(\SSS_{k,l}u)_k,z][(\SSS_{k,l}u)_l,z'] \\
	&=[(\SSS_{i,j}u)_i,z][(\SSS_{i,j}u)_j,z'] \\
	&=(\SSS^i(\SSS_{i,j}u)_i|z)(\SSS^j(\SSS_{i,j}u)_j|z').
\end{align*}
The other equality for case \ref{casoa} is analogous, using lemma \ref{SnSn}.
The same proof aplies to the case \ref{casob}, but replacing
$(\SSS_{i,j}(u))_i\otimes(\SSS_{i,j}(u))_j$ by
$(\SSS_{i,j}(u))_j\otimes(\SSS_{i,j}(u))_i$.
\end{pf}

We are in position now to give the last construction of a TOBA.

\begin{defn}
Let $V,W$ be $\k$-vector spaces with braidings $c_V,c_W$ and let
$(|):V\otimes W\to\k$ be a non degenerate bilinear form satisfying
\ref{casoa} (resp. \ref{casob}) of remark \ref{bfhe}. We take
$[,]:A(V)\otimes A(W)\to\k$ as in definition \ref{defcor}.
Let $\ide=\{v\in AV\ |\ [v,w]=0\ \forall w\in AW\}$,
$\ide'=\{w\in AW\ |\ [v,w]=0\ \forall v\in AV\}$ be the radicals of the
form $[,]$. Since $(|)$ is non degenerate, it is clear that
$\ide=\oplus_{n\ge 0}\ide_n=\oplus_{n\ge 0}\ker(\SSS^n:A^nV\to A^nV)$, and
$\ide'=\oplus_{n\ge 0}\ide'_n=\oplus_{n\ge 0}\ker(\SSS^n:A^nW\to A^nW)$.
Hence, another way to define $\toba V$ is to take $A(V)$ and divide out
by the left radical of the form $[,]$
(resp. for $\toba W$, we take $A(W)$ and divide out by the right radical).
\end{defn}

\begin{rem}
For the definition of the TOBA it is necessary for $(|)$ to be non degenerate,
though it is not necessary for the definition of the form $[,]$ in definition
\ref{defcor}.
\end{rem}

\begin{rem}
In the case of example \ref{ejelus} we get the algebra
$\f=\toba {\k(q)\theta_1\oplus\ldots\oplus\k(q)\theta_n}$.
For $V\in\yd$, $W=\dai V\in\yd$, we get the same object $\toba V$ as before,
and hence this is really an alternative form for the construction.
\end{rem}

Lemma \ref{relpc} says that the left and right radicals of $[,]$ are
Hopf ideals, and hence the relations between $AV$ and $AW$ give similar
relations between $\toba V$ and $\toba W$.

\begin{thm}
Let $V,W$ be as in lemma \ref{relpc}.
There is a unique non degenerate bilinear form
$\toba V\otimes\toba W\to\k$ such that
\begin{enumerate}
\item $[1,1]=1$, \label{punto1}
\item $[\tobag iV,\tobag jW]=0$ if $i\neq j$,
\item $[v,w]=(v|w)$ if $v\in\tobag 1V$, $w\in\tobag 1W$,
\item $[u,z\cdot z']=
	\left\{ \begin{array}{ll}
	[\com u1,z][\com u2,z'] & \mbox{ in case \ref{casoa}}, \\
	{}[\com u1,z'][\com u2,z] & \mbox{ in case \ref{casob}},
	\end{array}\right.$
\item $[u\cdot u',z]=
	\left\{ \begin{array}{ll}
	[u,\com z1][u',\com z2] & \mbox{ in case \ref{casoa}}, \\
	{}[u,\com z2][u',\com z1] & \mbox{ in case \ref{casob}},
	\end{array}\right.$ \label{punto5}
\end{enumerate}
where we denote $\tobag nV=A^nV/\ide_n$ the component in degree $n$.
\end{thm}

\begin{pf}
As in definition \ref{defcor},
we define the form $[,]=[,]_>$ in case \ref{casoa}, and $[,]=[,]_<$
in \ref{casob}. Lemma \ref{relpc} allows to consider
$[,]:\toba V\otimes\toba W\to\k$ induced from $[,]:AV\otimes AW\to\k$,
which turns out to be a non degenerate bilinear form satisfying
\ref{punto1}--\ref{punto5}. The uniqueness follows easily by induction.
\end{pf}

When $V,W$ are objects in $\cate$, a braided abelian category, and the
pairing $[,]:AV\otimes AW\to\k$ is a morphism in $\cate$, we have
a relation between $\toba W$ and $\dai{\toba V}$, provided this latter
object exists in $\cate$. We close the section giving the explicit
relation for $V$, $\dai V$ in $\yd$:

\begin{prop}
Let $V$ be an object in $\yd$. If $\toba V$ is finite dimensional, then
$\toba{\dai V}\simeq(\dai{\toba V})^{\hbox{\chq bop}}$.
\end{prop}

\begin{pf}
First, $\toba{\dai V}$ can be identified, via $[,]$, to $\dai{\toba V}$ as
a vector space. This identification is furthermore an isomorphism in $\yd$,
since $[,]$ is a morphism in $\yd$. We have to check the relation between
$[,]$ and the products and coproducts in $(\dai{\toba V})^{\hbox{\chq bop}}$,
$\toba V$ and $\toba{\dai V}$. We do it for the multiplication in
$\toba{\dai V}$, the other one being analogous.

Let $\{\bdr\}_{\alpha}$, $\{\dbr\}_{\alpha}$ be dual bases of $\toba V$.
We have for $u\in\toba V$, $f,g\in(\dai{\toba V})^{\hbox{\chq bop}}$,
we have
\begin{align*}
\gi u,&m_{\hbox{\chq bop}}(f\otimes g)\gd \\
&=\gi u,(\con f1g)\com f0\gd \\
&=\gi\com u1,\comn u21\con f1g\gd\gi\comm u20,\com f0\gd \\
&=\gi\com u1,\comn u21\Ss(\con{\bdr}1)g\gd \\
	&\qquad\gi\comm u20,\dbr\gd\gi\com{\bdr}0,f\gd \\
&=\gi\com u1,\comn u22\Ss(\comn u21)g\gd\gi\comm u20,f\gd \\
&=\gi\com u1,g\gd\gi\com u2,f\gd,
\end{align*}
but this is the same equality for $\toba{\dai V}$.
\end{pf}

\subsection{Concrete examples}\ \\
We present now two families of braided Hopf algebras discovered by
Milinski and Schneider. Both families are
particular cases of Hopf algebras in Yetter--Drinfeld categories over
group algebras of Coxeter groups (see \cite{milinski-schneider})
and have the form $A(V)/I$ for certain $V$ and $I\subset\ker(\SSS)$.
It is not known in general whether or not $I=\ker(\SSS)$ (that is, whether
or not they are TOBAs). Most of the results in this section are taken from
to \cite{milinski-schneider}, an exception is Proposition \ref{infdim}.

\begin{exmp}
Let $n\in\NN$, and $H=\k\SS_n$. We take $V$ the $\k$-vector space
with basis consisting of elements $y_{\tau}$ where $\tau$ runs over all
(not only elementary) transpositions $\tau=(i,j),\ i\neq j$.
We make $V$ an object of $\yd$ taking
\begin{align*}
\delta(y_{\tau})&=\tau\otimes y_{\tau}, \\
\sigma\ganda y_{\tau}&=\sg(\sigma)y_{\sigma\tau\sigma^{-1}}.
\end{align*}
The module $V$ is nothing but $M(g,\rho)$ with $g$ any transposition and
$\rho$ the restriction of the sign representation to the isotropy
subgroup of $g$.

Let now $J$ be the subspace of $A^2(V)$ generated by the elements
\begin{align}
\label{sn1} y_{\tau}^2 \qquad &\forall\tau,\\
\label{sn2} y_{\tau}y_{\tau'}+y_{\tau'}y_{\tau} \qquad
		&\mbox{if }\tau\tau'=\tau'\tau, \\
\label{sn3} y_{\tau}y_{\tau'}+y_{\tau'}y_{\tau''}+y_{\tau''}y_{\tau} \qquad
		&\mbox{if } \tau\tau'=\tau''\tau.
\end{align}
Then $J=\ker(\SSS^2)$. Let $I$ be the ideal generated by $J$.
Since $J$ is an $H$-submodule and an $H$-subcomodule, the same
is true for $I$. Since $J$ is a coideal, the same is true for $I$.
Then $R^1_n:=A(V)/I$ is a Hopf algebra in $\yd$.
\end{exmp}

\begin{exmp}
Let $p$ be an odd prime number. We take now $H=\k\DD_p$, where
$\DD_p$ is the dihedral group, i.e. the group generated by $\rho$ and
$\sigma$, with relations
$$\rho^p=\sigma^2=1,\quad \sigma\rho=\rho^{p-1}\sigma.$$
The conjugacy class of $\sigma$ is
$\orb_\sigma=\{\sigma,\rho\sigma,\ldots,\rho^{p-1}\sigma\}$,
and the isotropy subgroup is $(\DD_p)_\sigma=\{1,\sigma\}$. We take now
$\chi:(\DD_p)_\sigma\to\k^\times,\ \chi(\sigma)=-1$,
and then define $V=M(\sigma,\chi)$ as in \ref{loseme}. 
Let $V_0$ be the space affording $\chi$. We denote $y_0$ a generator of $V_0$.
We put $y_i=\rho^i\ganda y_0\in M(\sigma,\chi)$. We take the subindices of
the $y_i$ to be on $\ZZ/p$, thus $y_{i+p}=y_i$. We then have
\begin{align*}
&\delta(y_i)=\alai \sigma{\rho^i}\otimes y_i=\rho^{2i}\sigma\otimes y_i, \\
&\rho^j\ganda y_i=y_{i+j},\qquad \sigma\ganda y_i=-y_{-i}, \\
&c_{V,V}(y_i\oti y_j)=-y_{2i-j}\oti y_i.
\end{align*}
To compute $\ker\SSS^2$ it is convenient to take a different basis in
$V\otimes V$. Let $\xi$ be a primitive $p$-root of unit (we may suppose
$\k$ has a primitive $p$-root of unit, for, if not, we can take a
suitable extension of $\k$). Let
$$w^r_k=\sum_{i=0}^{p-1}\xi^{ri}y_i\oti y_{i+k},\qquad
	0\le r,k<p.$$
Then the $w^r_k$ form a basis of $V^{\otimes 2}$, and
$c_{V,V}(w^r_k)=-\xi^{rk}w^r_k$, whence
\begin{equation} \label{kers2}
\ker(\SSS^2)=\ker(\id+c)=\gi w^r_k,\ rk=0\gd
	\qquad\mbox{(and then $\dim\ker(\SSS^2)=2p-1$)}.
\end{equation}
It is easy to see that $\gi w^0_0,w^1_0,\ldots,w^{p-1}_0 \gd=
\gi (y_0\oti y_0),\ (y_1\oti y_1),\ldots,\ (y_{p-1}\oti y_{p-1})\gd$.
We define then $R^2_p$ to be $A(V)/I$, where $I$ is generated by
$$\begin{array}{rl}
& y_i\oti y_i,\qquad 0\le i<p, \\
w^0_1=&y_0\oti y_1 + y_1\oti y_2 + \cdots + y_{p-1}\oti y_0, \\
w^0_2=&y_0\oti y_2 + y_1\oti y_3 + \cdots + y_{p-1}\oti y_1, \\
\hdotsfor{2} \\
w^0_{p-1}=&y_0\oti y_{p-1} + y_1\oti y_0 + \cdots + y_{p-1}\oti y_{p-2}.
\end{array}$$
\end{exmp}

\begin{rem}
The hypothesis $p$ being an odd prime number is not necessary. It is
used because it makes the relations simpler.
\end{rem}

The algebras $R^2_p$ are infinite dimensional for $p>7$. This is a
consequence of the following
\begin{thm}[Golod--Shafarevich]
Let $V=\oplus_{n>0}V_n$ be a graded vector space, and $A=T(V)$ be the
(graded) tensor algebra of $V$. Let $I$ be a homogeneous ideal, and suppose
$I$ is generated (as an ideal) by the subspaces $\oplus_{n>0}I_n$.
Let $R=T(V)/I$ be the quotient algebra. Let $h_V$ and $h_I$ be the
Hilbert series of $V$ and $I$, that is,
$$h_V(t)=\sum_{n>0}\dim(V_n)t^n,\qquad h_I(t)=\sum_{n>0}\dim(I_n)t^n.$$
Let $g(t)=\sum g_nt^n=(1-h_V(t)+h_I(t))^{-1}$ as formal power series.
If $g_n\ge 0\ \ \forall n$, then $R$ is infinite dimensional.
\end{thm}

\begin{pf} See \cite{ufnarovski}. \end{pf}

\begin{prop}\label{infdim}
The algebras $R^2_p$ are infinite dimensional for $p>7$.
\end{prop}

\begin{pf}
We apply the theorem. We have $h_V(t)=pt$, and $h_I(t)=(2p-1)t^2$.
Then
$$g(t)=(1-pt+(2p-1)t^2)^{-1}=(1-t/a)^{-1}(1-t/b)^{-1}=
	(\sum_{n\ge 0}(t/a)^n)(\sum_{n\ge 0}(t/b)^n)$$
for $a,b$ the roots of $(1-pt+(2p-1)t^2)$. If $a$ and $b$ are both real
and positive then $g_n\ge 0\ \forall n$. This is true if $p^2-4(2p-1)\ge 0$,
which implies $p>7$.
\end{pf}

\bigskip

For $p=3$ we have $\DD_3\simeq\SS_3$, and in fact $R^2_3\simeq R^1_3$.

\begin{prop} \label{erretres}
$R^2_3\simeq R^1_3$ is a TOBA of dimension $12$.
\end{prop}

\begin{pf}
It can be seen by direct computation using the relations that
\begin{align*}
y_0y_1y_0 &= -y_1y_2y_0 = y_1y_0y_1 = -y_0y_2y_1, \\
y_0y_1y_2 &= -y_0y_2y_0 = y_2y_1y_0 = -y_2y_0y_2, \\
y_1y_0y_2 &= -y_2y_1y_2 = y_2y_0y_1 = -y_1y_2y_1,
\end{align*}
and the other monomials in degree $3$ vanish since in all of them
appears $y_i^2$ for some $i$. This in turn implies
\begin{equation} \label{relag4}
\begin{split}
y_0y_1y_0y_2 &= -y_1y_2y_0y_2 = y_1y_0y_1y_2 = -y_0y_2y_1y_2 =
	y_0y_2y_0y_1 = -y_0y_1y_2y_1 \\
&= -y_2y_1y_0y_1 = y_2y_0y_2y_1 = -y_2y_0y_1y_0 = -y_1y_0y_2y_0 =
	y_2y_1y_2y_0 = y_1y_2y_1y_0,
\end{split}
\end{equation}
\begin{align*}
y_0y_1y_0y_1 &= y_1y_2y_0y_1 = y_1y_0y_1y_0 = y_0y_2y_1y_0 =
	y_0y_1y_2y_0 = y_2y_0y_2y_0 = y_0y_2y_0y_2 \\
&= y_2y_1y_0y_2 = y_1y_0y_2y_1 = y_2y_1y_2y_1 = y_2y_0y_1y_2 =
	y_1y_2y_1y_2 = 0,
\end{align*}
and the other monomials in degree $4$ vanish since in all of them
appears $y_i^2$ for some $i$.
Moreover, the monomials in \eqref{relag4} are annihilated by multiplying
them with any of the $y_i$, and then
$$\agg{R^2_3}n=0\ \forall n\ge 5.$$
With this, we get the set of generators of $R^2_3$ consisting of
\begin{equation} \label{congr3}
\{1,\ y_0,\ y_1,\ y_2,\ y_0y_1,\ y_1y_2,\ y_0y_2,
\ y_1y_0,\ y_0y_1y_0,\ y_0y_1y_2,\ y_1y_0y_2,\ y_0y_1y_0y_2\}.
\end{equation}
It can be proved that this set is indeed a basis taking the representation
of rank $12$ given by
\begin{align*}
y_0 &\mapsto A_0 = E_{1,2}+E_{3,7}+E_{4,8}+E_{5,9}+E_{6,10}+E_{11,12}; \\
y_1 &\mapsto A_1 = E_{1,3}+E_{2,5}-E_{4,6}-E_{4,7}-E_{6,9}+E_{7,9}-
	E_{8,11}+E_{10,12}; \\
y_2 &\mapsto A_2 = E_{1,4}+E_{2,6}-E_{3,5}-E_{3,8}-E_{5,10}-E_{7,11}+
	E_{8,10}+E_{9,12};
\end{align*}
where $E_{i,j}$ stands for the matrix with $1$ in the entry $(i,j)$ and
$0$ in the others. This is easily seen to be a representation
(i.e. $A_0^2=A_1^2=A_2^2=A_0A_1+A_1A_2+A_2A_0=A_0A_2+A_1A_0+A_2A_1=0$)
and to map the set in \eqref{congr3} to a linearly independent set,
which says that $\dim R^2_3=12$. Alternatively one can use the
Diamond Lemma.

\medskip
We have to check now that $R^2_3$ is a TOBA. Let $V=\agg{R^2_3}1$, and
let $T=\toba V$. Since $I\subseteq\ker\SSS$, we know that there exists
a surjective graded morphism $\pi:R^2_3\sobre T$. Let $N$ be such that
$\agg TN\neq 0$ and $\agg Ti=0\ \forall i>N$. By \ref{dualpoin}
we have that $\dim\agg TN=1$ and $\dim\agg Ti=\dim\agg T{N-i}$.
We have then the following possibilities:
\begin{enumerate}
\item $N=4$, and then $\dim\agg T3=\dim\agg T1=\dim V=3$, from where $\pi$
is an isomorphism unless $\dim\agg T2<4$.
\item $N=3$, and then $\dim\agg T2=\dim\agg T1=3$.
\item $N=2$, and then $\dim\agg T2=\dim\agg T0=1$.
\end{enumerate}
We see that in any case $\pi$ is an isomorphism unless $\dim\agg T2<4$,
but $\dim\agg T2$ is the codimension of $\ker\SSS^2$ in $V\otimes V$,
and we know from \ref{kers2} that it is equal to $4$.
\end{pf}

\medskip
We give now the bosonised algebra. We denote also by $y_i$ the element
$(y_i \# 1)$, and by $g_0,\ g_1$ the group-likes
$(1\# \sigma),\ (1\# \rho^2\sigma)$ (which generate $\DD_3$).
The bosonization is thus the algebra presented by generators
$g_0,\ g_1,\ y_0,\ y_1,\ y_2$ with relations
\begin{align}
\label{rt1} & g_{i}^2 = 1 \qquad i = 0,1;\\
\label{rt2} & g_{1}g_{0}g_{1} = g_{0}g_{1}g_{0}; \\
\label{rt3} & y_{j}^2 = 0 \qquad j = 0,1,2;\\
\label{rt4} & y_{0}y_{1} + y_{1}y_{2} + y_{2}y_{0} = 0; \\
\label{rt5} & y_{1}y_{0} + y_{0}y_{2} + y_{2}y_{1} = 0;\\
\label{rt6} & g_{0}y_{0}g_{0} = -y_{0}, \quad g_{0}y_{1}g_{0} = -y_{2},
	\quad g_0y_2g_0 = -y_1; \\
\label{rt7} & g_{1}y_{0}g_{1} = -y_{2}, \quad g_{1}y_{1}g_{1} = -y_{1},
	\quad g_1y_2g_1 = -y_0.
\end{align}
The Hopf algebra structure is determined by 
$$\Delta(g_{i}) = g_{i} \otimes g_{i} \quad (i=0,1),
	\qquad \Delta(y_i) = y_i\otimes 1+g_i\otimes y_i \quad (i=0,1,2),$$
where we denote $g_2=g_0g_1g_0$. This Hopf algebra has dimension 72; it is
pointed and its coradical is isomorphic to the group algebra of $\SS_3$.

\bigskip
\begin{rem}
More Hopf algebras with coradical $\k\SS_3$ appear replacing the relations
\eqref{rt4} and \eqref{rt5} by
\begin{align*}
y_{0}y_{1} + y_{1}y_{2} + y_{2}y_{0} &= \lambda_1(g_0g_1-1), \\
y_{1}y_{0} + y_{0}y_{2} + y_{2}y_{1} &= \lambda_2(g_1g_0-1),
\end{align*}
with $\lambda_1,\lambda_2\in\k$. We will consider this and related problems
in a separated article.
\end{rem}

\bigskip
\begin{rem}
The relations \eqref{rt3} can be twisted as in the preceding remark,
but in this case one must replace the group $\SS_3$ by a covering of it,
using the relations $g_i^{2N}=1$ instead of $g_i^2=1$.
\end{rem}

\bigskip
It is shown in \cite{milinski-schneider} that $R^1_4$ is finite dimensional.
It is not known whether $R^1_n$ is finite dimensional or not for $n>4$.

\bigskip
It is not known whether $R^2_5$ and $R^2_7$ are finite dimensional or not.
It is not known neither whether the algebras obtained dividing $A(V)$ by
$\ker(\SSS)$ (and not only by the ideal generated by $\ker(\SSS^2)$) are finite
dimensional or not.

\bigskip
\begin{exmp}
As a last example, we take $\Gamma=\DD_4$, $H=\k\Gamma$. The conjugacy
class of $\sigma$ is $\orb_\sigma=\{\sigma,\rho^2\sigma\}$, and the conjugacy
class of $\rho\sigma$ is $\orb_{\rho\sigma}=\{\rho\sigma,\rho^3\sigma\}$.
We take then, in a similar way to $R^2_p$, the module $V$ in $\yd$ with
basis $\{z_0,z_1,z_2,z_3\}$ with the structure given by
$\delta(z_i)=\rho^i\sigma\otimes z_i$, $\rho^j\ganda z_i=z_{i+2j}$ and
$\sigma\ganda z_i=-z_{-i}$ (where as before we take the subindices of the
$z_i$ to be on $\ZZ/4$). Then $V$ can be decomposed as
$$V=\gi z_0,z_2\gd \oplus \gi z_1,z_3\gd = V_0\oplus V_1,$$
which are irreducible. We have as before that the braiding is given by
$$c(z_i\otimes z_j)=-z_{2i-j}\otimes z_i.$$
Let $T_0=\toba{V_0}$ and $T_1=\toba{V_1}$. It is easy to see that
$T_i\ (i=0,1)$ have dimension $4$, and their respective ideals $\ker\SSS$
are generated by
\begin{align*}
z_0^2=z_2^2=0, \quad & z_0z_2+z_2z_0=0, \\
z_1^2=z_3^2=0, \quad & z_1z_3+z_3z_1=0.
\end{align*}
The TOBA $T=\toba V$ is much more complicated to compute.
Let $a=z_1z_2+z_0z_1$ and $b=z_1z_0+z_2z_1$.
Then the elements $a^2,\ b^2$ and $ab+ba$ are primitive in $A(V)$,
which says that they belong to $\ker\SSS$. We have then a graded
braided Hopf algebra dividing out $A(V)$ by the relations
\begin{align*}
& z_0^2=z_2^2=0, \quad z_0z_2+z_2z_0=0, \\
& z_1^2=z_3^2=0, \quad z_1z_3+z_3z_1=0, \\
& z_0z_1+z_1z_2+z_2z_3+z_3z_0=0, \\
& z_0z_3+z_1z_0+z_2z_1+z_3z_2=0, \\
& a^2=b^2=0,\quad  ab+ba=0.
\end{align*}
Using the Diamond Lemma, it can be shown that the dimension of this algebra
is $64$. This is done in \cite{milinski-schneider}.
\end{exmp}

\end{document}